\documentclass[11pt]{amsart}
\usepackage[UKenglish]{isodate}
\usepackage{amsmath}
\usepackage{amsthm}
\usepackage{amsbsy}
\usepackage{amssymb}
\usepackage{amsfonts}
\usepackage{appendix}
\usepackage[dvips]{lscape}
\usepackage{xcolor}
\usepackage{amscd}
\usepackage[all,cmtip]{xy}
\usepackage{euscript}
\usepackage{caption}  
\usepackage{parskip}
\usepackage{enumerate}
\usepackage{yfonts}
\usepackage{xcolor}
\usepackage{graphicx}
\usepackage{fancyhdr}
\usepackage{tikz-cd}
\usepackage{colortbl}
\usepackage{comment}
\usepackage{fullpage}
\usepackage{verbatim}
\usepackage{subcaption}

\usetikzlibrary{positioning}
\usetikzlibrary{decorations.pathreplacing}

\usepackage[margin=1in]{geometry}
\usepackage[expansion=false]{microtype}

\theoremstyle{plain}

\theoremstyle{plain}
\newtheorem{theorem}{Theorem}[section]

\newtheorem{proposition}[theorem]{Proposition}
\newtheorem{lemma}[theorem]{Lemma}
\newtheorem{corollary}[theorem]{Corollary}

\newcommand{\res}{\text{Res}}

\theoremstyle{remark}
\newtheorem{remark}[equation]{Remark}

\theoremstyle{definition}

\newcommand{\vol}{\textup{vol}}

\newif\iffinalrun
\iffinalrun
\else
 \fi

\iffinalrun
  \newcommand{\need}[1]{}
  \newcommand{\mar}[1]{}
\else
  \newcommand{\need}[1]{{\tiny *** #1}}
  \newcommand{\mar}[1]{\marginpar{\raggedright\tiny Fix Me:  #1 }}\fi

\usepackage{amssymb,amsmath,amsfonts,amsthm,epsfig,amscd,stmaryrd}
\usepackage{stmaryrd}

\usetikzlibrary{arrows}

\pgfdeclarelayer{background}
\pgfsetlayers{background,main}

\pgfmathsetmacro{\myxlow}{-2}
\pgfmathsetmacro{\myxhigh}{2}
\pgfmathsetmacro{\myiterations}{6}

\title{ New Lower Bounds for the Schur-Siegel-Smyth Trace Problem }
\author{Bryce Joseph  Orloski, Naser Talebizadeh Sardari and Alexander Smith}

\address{Penn State department of Mathematics, McAllister Building, Pollock Rd, State College, PA 16802 USA}
\email{nzt5208@psu.edu,bjo5149@psu.edu}

\address{UCLA Department of Mathematics, 520 Portola Plaza, Los Angeles, CA 90095 USA}
\email{asmith13@math.ucla.edu}
\thanks{We would like to thank Professors  Sarnak, Zargar, and Zikatanov for their comments on the earlier version of this work. The authors also thank Pennsylvania State University's Institute for Computational and Data Sciences for allowing us to use their ROAR servers for our numerical experiments.
This work was partially supported by NSF grant DMS-2401242.
 The third author served as as Clay Research
Fellow during the writing of this paper, and would like to thank the Clay Mathematics Institute
for their support.
}

\begin{document}
\begin{abstract}
    We derive and implement a new way to find lower bounds on the smallest limiting trace-to-degree ratio of totally positive algebraic integers and improve the previously best known bound to 1.80203. Our method adds new constraints to Smyth's linear programming method to decrease the number of variables required in the new problem of interest. This allows for faster convergence recovering Schur's bound in the simplest case and Siegel's bound in the second simplest case of our new family of bounds. 
    We also prove the existence of a unique optimal solution to our newly phrased problem and express the optimal solution in terms of polynomials. 
    Lastly, we solve this new problem numerically with a gradient descent algorithm to attain the new bound 1.80203.
\end{abstract}

\maketitle

\section{Introduction}\label{intro}

\subsection{The Schur-Siegel-Smyth trace problem.}
For an algebraic integer $\alpha$ with $\deg(\alpha)=n$ and the set of Galois conjugates   $[\alpha]:=\{\alpha_1,\dots,\alpha_n\}$, define 
\[\overline{\text{tr}}(\alpha):=\frac{\sum_{i=1}^{n}\alpha_i}{n},\]
\[ \text{Nr}(\alpha)=\prod_{i=1}^n\alpha_i,\text { and   }\]
\[ \Delta(\alpha)=\prod_{i<j} (\alpha_i-\alpha_j)^2\in \mathbb{Z}^+. \]

The algebraic integer $\alpha$ is called totally positive if $\alpha_i>0$ for $i=1,\dots, n$.
Let $\mathcal A$ be the set of totally positive algebraic integers.
The Schur-Siegel-Smyth trace problem is to compute the value of the number $\lambda^{SSS}:=\liminf_{\alpha\in\mathcal A}\overline{\text{tr}}(\alpha)$ for a fixed ordering of $\mathcal A$.
Schur~\cite{Schur} in 1918 proved that
\[
\lambda^{SSS}\geq \sqrt{e}.
\] 
In other words, this means that for any $\lambda < \sqrt{e}$, there are only finitely many $\alpha\in\mathcal A$ for which $\overline{\text{tr}}(\alpha) < \lambda$.
Siegel~\cite{MR12092} improved this  to
\[
\lambda^{SSS}\geq e(1+\nu^{-1})^{-\nu}
\approx 1.7336105,
\]
where
$\nu$ is the positive root of the transcendental equation
\begin{equation*}
(1+\nu)\log(1+\nu^{-1})+\frac{\log \nu}{1+\nu}=1.
\end{equation*}
Siegel used only that $\text{Nr}(\alpha)\geq 1$ and $\Delta(\alpha)\geq 1$ for $\alpha\in\mathcal A$, and he solved a similar optimization problem as Schur using Lagrange multipliers.
\newline

The previous best lower bound was achieved numerically by using a linear programming method initiated first by Smyth~\cite{MR0736460}. In what follows, we explain Smyth's method. 
Smyth showed that if $\mu$ is the limit of uniform probability measures attached to distinct conjugate algebraic integers 
all of whose conjugates lie in a compact set $K$ of the complex plane, then 
\[
 \int \log|Q(x)|d\mu(x)\geq 0
\]
for every $Q(x)\in \mathbb{Z}[x]$.
In a recent breakthrough~\cite{Smith}, the third author proved that these necessary conditions are sufficient when $K$ is 
a subset of $\mathbb{R}$ under some assumptions on $K$. The first two authors recently established this result in \cite{OT} for compact subsets $K$ of the complex plane under some assumptions, and the new method gives a different treatment of Smith's Theorem.
In particular, we introduced some new constraints on the limiting measures of conjugate algebraic integers.  For a limiting  measure $\mu$, the first two authors proved in~\cite{OT} that
\begin{equation}\label{multi}
\int \log(|Q(x_1,\dots,x_n)|)d\mu(x_1)\dots d\mu(x_n)\geq 0
\end{equation}
for every $Q(x_1,\dots,x_n)\in \mathbb{Z}[x_1,\dots,x_n]$. The goal of this paper is to improve Smyth's lower bound method by using new constraints in~\eqref{multi}.

\subsection{Main Results}\label{prime}
Smyth gave a lower bound on the trace  problem by minimizing $\int x d\mu(x)$ subject to 
\[\int \log|Q(x)|d\mu(x)\geq 0\]
for every $Q(x)\in A$, where $A\subset\mathbb{Z}[x] $ is  a fixed finite subset.
In fact, Smyth worked with the equivalent associated dual linear programming problem~\cite[Chapter 5]{Boyd} by finding maximal $\lambda_A$ such that
\begin{equation}\label{1var}
x\geq \lambda_A+\sum_{Q\in A}\lambda_Q \log|Q(x)|    
\end{equation}
holds for every $x\in \mathbb{R}^+$ and some choice of scalars $\lambda_Q\geq 0$.  In this paper, we minimize $\int x d\mu(x)$ subject to inequalities \eqref{multi} for every $Q\in \mathbb{Z}[x_1,\dots,x_n]$. Following Smyth's method, we consider maximizing   $\lambda_B$ such that
\begin{equation}\label{2var}
\frac{x+y}{2}\geq \lambda_B+\sum_{Q\in B} \lambda_Q \log|Q(x,y)|    
\end{equation}
for every $x,y\in\mathbb R^{+}$ and some $\lambda_B,\lambda_Q\geq 0$ and a finite subset $B\subset \mathbb{Z}[x,y]$.
By integrating both sides of~\eqref{2var} with respect to $d\mu(x)d\mu(y)$, it follows that
\[
\int x d\mu(x)\geq \lambda_B.
\]
We note that the bounds in \eqref{2var} are  at least as good as Smyth's method which uses~\eqref{1var}. We verify this next. Suppose that~\eqref{1var} holds for some choices of $A\subset \mathbb{Z}[x]$ and $\lambda_A,\lambda_Q\geq 0$. This implies
\[
\begin{cases}
    x\geq \lambda_A+\sum_{Q\in A}\lambda_Q \log|Q(x)|,\\
     y\geq \lambda_A+\sum_{Q\in A}\lambda_Q \log|Q(y)|.
\end{cases}
\]
By taking the mean of the above inequalities, we have 
\[
\frac{x+y}{2}\geq \lambda_A+\sum_{Q\in A}\lambda_Q \frac{\log|Q(x)|+\log|Q(y)|}{2} 
\]
for every $x,y\in\mathbb R^{+}$. This is a special case of \eqref{2var} with
\[
B:=\{Q(x),Q(y): Q\in A \}.
\]
This implies that $\lambda_A\leq \lambda_B$. Moreover, 
our main theorem shows that Smyth's bound $\lambda_A$ is strictly smaller than $\lambda_{B^+}$, where $B^+=B\cup\{x-y\}$. 
\\

There is a loss in working with the dual problem~\eqref{2var} as the constraints in~\eqref{multi} are not linear in $\mu$. Next, we improve upon this dual bound by working directly with~\eqref{multi}.
The main difficulty in working directly with~\eqref{multi} is that, as stated before, this problem is no longer a linear programming problem and the associated dual problem in \eqref{2var} is not equivalent in general. Our main idea is to find a particular $Q(x_1,\dots,x_n)\in \mathbb{Z}[x_1,\dots,x_n]$ such that the functional $\mu \mapsto \int \log(|Q(x_1,\dots,x_n)|)d\mu(x_1)\dots d\mu(x_n)$ is concave, and then apply the strong duality identity which holds for the convex optimization problems. 
We state  our main theorem  for $Q(x_1,x_2)=x_1-x_2$. A similar method works for $Q(x_1,\dots,x_4)=x_1-x_2-x_3+x_4$ which is a possible object of future research.
In the case of $Q(x_1,x_2)=x_1-x_2$, we define
\[
I(\mu):=\int \log|x_1-x_2|d\mu(x_1)d\mu(x_2)
\]
to be the logarithmic energy of $\mu$.

\begin{theorem}\label{main1} 
    Suppose that $A\subset \mathbb{R}[x]$ is a finite subset of polynomials with real roots. Let $\lambda_A= \inf_{\mu} \int x d\mu(x)$, where $\mu$ varies  among probability measures supported on $\mathbb{R}^+$ with $I(\mu)\geq 0$ and $\int \log|Q(x)|d\mu \geq 0$ for every $Q\in A$. There exists a unique measure $\mu_A$ supported on a finite union of intervals $\Sigma\subset [0,18]$ such that
    \[
    \begin{array}{cc}
          \int x d\mu_A(x)=\lambda_{A},\\
         I(\mu_A)=0,          \\
                   x\geq \lambda_A +\sum_{Q\in A}\lambda_{Q} \log|Q(x)|+\lambda_0 U_{\mu_A}(x),\text{ and}\\
          \int \log|Q(x)| d \mu_A=0 \text{ if $\lambda_Q\neq 0$,}

    \end{array}
    \]
where in the third line the inequality only holds for non-negative $x$ and equality holds for every $x\in \Sigma$ with some scalars $\lambda_Q\geq 0$ and $\lambda_0> 0$. Moreover if $\Sigma:=\bigcup_{i=0}^l[a_{2i},a_{2i+1}] $ and $\{x\}\subset A$, then $\mu_A$ has the following density function up to a scalar
\[
\frac{|p(x)|\sqrt{|H(x)|}}{\prod_{\lambda_Q\neq 0}|Q(x)|},
\]
where $H(x)=\prod_{i=0}^l(x-a_{2i})(x-a_{2i+1})$ and $p(x)$ is a polynomial with degree $\sum_{\lambda_Q\neq 0} \deg(Q)-l-1$. Furthermore, every gap interval $(a_{2i+1},a_{2i+2})$ contains a root of some $Q\in A$.
\end{theorem}
\begin{corollary}
    Let $\lambda_A$ be as in Theorem~\ref{main1}. We have
    \[
    \lambda_A\le \lambda^{SSS}.
    \]
\end{corollary}
\begin{proof}
    Suppose that $\alpha_n$ is a sequence of totally positive algebraic integers such that
    \[
    \lim_{n\to \infty}\overline{\text{tr}}(\alpha_n)=\lambda^{SSS}.
    \]
    Let 
    \[
    \mu_{\alpha_n}:=\frac{1}{\deg(\alpha_n)}\sum_{\alpha_{n,j}\in [\alpha_n]}\delta_{\alpha_{n,j}}
    \]
    be the uniform discrete probability measure on the set of Galois conjugates  of $\alpha_n$. Since $\mu_{\alpha_n}$ is a tight family of probability measures, Prokhorov's theorem implies that there exists a convergent subsequence. Suppose that
\(
\mu
\)
is a limiting measure of a subsequence. 
As we mentioned before, $\mu$ satisfies all the constraints of Theorem~\ref{main1} for any finite set $A$. Hence,
\[
\lambda_A\leq \lambda^{SSS}
\]
for any finite set $A$.
\end{proof}
\begin{remark}
    Note that Corollary 5.5 of the third author in \cite{Smith} directly implies that the least upper bound of such $\lambda_{A}$ is $\lambda^{SSS}$.
\end{remark}

We recover Schur's lower bound by taking $A=\emptyset$ \cite{Schur} and Siegel's lower bound  by taking $A=\{x\}$ \cite{MR12092}. We state this in the following corollaries.
\begin{corollary}\label{schur_bound}
    If $A=\emptyset$, then $\lambda_A = \sqrt e$. Moreover, $d\mu_A = \frac{1}{2\pi}\sqrt{\frac{(4\sqrt{e}-x)}{ex}}dx$.
\end{corollary}

\begin{corollary}\label{siegel_bound} If $A=\{x\}$, then
\[
\lambda_A=e(1+\nu^{-1})^{-\nu}=1.7336105\dots,
\]
where $\nu$ is the positive root of the transcendental equation
\begin{equation}\label{siegel_bound_eq}
(1+\nu)\log(1+\nu^{-1})+\frac{\log \nu}{1+\nu}=1.
\end{equation}
Moreover, 
\[
\mu_A=\frac{2\sqrt{(b-x)(x-a)}}{\pi(a+b-2\sqrt{ab})x}dx
\]
where
\[
a=e(1+\nu^{-1})^{-\nu}\left(1-\frac{1}{\sqrt{\nu+1}}\right)^2\text{ and } 
b=e(1+\nu^{-1})^{-\nu}\left(1+\frac{1}{\sqrt{\nu+1}}\right)^2.
\]
\end{corollary}

\subsection{Gradient Descent}\label{grad_desc}
We developed an algorithm which given a finite $A\subset \mathbb Z[x]$, returns a distribution approximating $\mu_A$.
Our algorithm is based on a gradient descent algorithm.

\begin{corollary}\label{2poly_bound}
Let $A=\{x,x-1\}$. Then
\[
|\lambda_A-1.7773797|\leq 10^{-7}.
\]
Moreover, $\mu$ is supported on the union of two intervals $[a_0,a_1]\cup[a_2,a_3]$ with the following numerical values:
\[
a_{0}=0.0362736\dots, a_{1}=0.828301\dots,a_{2}=1.190973\dots,
\text{ and }a_{3}=5.707091\dots.
\]
Furthermore, it has the following density function:
\[
d\mu_A(x)=\frac{\sqrt{|\left(x-a_{0}\right)\left(x-a_{1}\right)\left(x-a_{2}\right)\left(x-a_{3}\right)|}}{c|x\left(x-1\right)|}dx
\]
where
\[
c=6.420592\dots.
\]
\end{corollary}

\begin{corollary}\label{3poly_bound}
If $A = \{x, x-1, x^2-3x+1\}$, then
\[
|\lambda_A-1.793023|\leq 10^{-6}
\]
where $\Sigma$ is a union of four intervals approximated by
\[
  [0.0409275, 0.34114487]\cup [0.4252603, 0.811681]\cup [1.211488, 2.4844644] \cup [2.7580631, 5.52512172].
\]  
\end{corollary}

\begin{corollary}\label{4poly_bound}
    If $A=\{x,x-1,x^2-3x+1,x^3-5x^2+6x-1\}$, then 
    \[|\lambda_A- 1.798249|\leq 10^{-6}.\] 
\end{corollary}
\begin{corollary}\label{5poly_bound}
    If $A=\{x,x-1,x-2,x^2-3x+1,x^3-5x^2+6x-1\}$, then 
    \[|\lambda_A- 1.7998|\leq  10^{-4}.\] 
\end{corollary}
\begin{corollary}\label{our_bound}
    If $A=\{x,x-1,x-2,x^2-3x+1,x^3-5x^2+6x-1, x^4-7x^3+13x^2-7x+1,$ $x^4-7x^3+14x^2-8x+1\}$, then 
    \[|\lambda_A- 1.80208|\leq  5\times10^{-5}.\] 
\end{corollary}
 \begin{figure}
         \centering
    \includegraphics[scale=0.86]{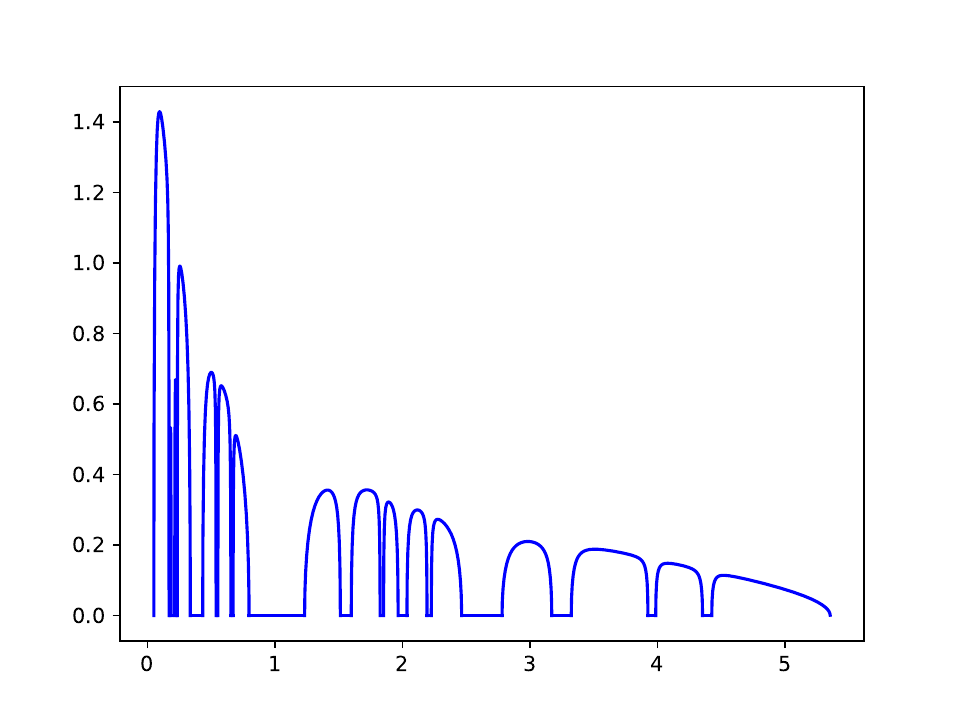}
         \caption{This is the density function of $\mu_{A}$ where $A$ is defined in Corollary \ref{our_bound}. Note that the two darker vertical lines near 0 are densities intervals of width too small to be accurately depicted.}
        \label{new-l_density}
     \end{figure}
This beats the currently best known lower bound on the problem. 
For extra details about the measure used in this corollary, see Figure~\ref{new-l_density} and Appendix \ref{numerical_measures}.
\newline

Below is a table of the history of the best known lower bound on the problem with the last line including the newest bound which is the main topic of this paper. This is the largest improvement on the lower bound since Smyth~\cite{MR0736460} in 1984.

\begin{center}
\begin{tabular}{|c|c|c|}
\hline Authors & Year & Bound Estimate\\\hline\hline
Schur~\cite{Schur}& 1918 & 1.6487\\\hline
Siegel~\cite{MR12092}& 1945 & 1.7336\\\hline
Smyth~\cite{MR0736460}& 1984 & 1.7719\\\hline
Flammang, et al.~\cite{flammang97} & 1997 & 1.7735\\\hline
McKee and Smyth~\cite{mckee_smith04}& 2004 & 1.7783\\\hline
Aguirre, Bilbao, and Peral~\cite{aguirre06} &  2006 & 1.78002\\\hline
Aguirre and Peral~\cite{aguirre07} & 2007 & 1.7836\\\hline
Aguirre and Peral~\cite{MR2428512}& 2008 & 1.7841\\\hline
Flammang~\cite{flammang09}& 2009 & 1.78702\\\hline
McKee~\cite{McKee}& 2011 & 1.78839\\\hline
Liang and Wu~\cite{Liang} & 2011 & 1.7919\\\hline
Flammang~\cite{2019arXiv190709407F} & 2016 & 1.7928\\\hline
Wang, Wu, and Wu~\cite{MR4280303} & 2021 & 1.7931\\\hline
Our new lower bound & 2023 & 1.80203    \\\hline
\end{tabular}
\newline
\end{center}

The proof of the previously best known lower bound by Wang, Wu, and Wu uses Smyth's method with $130$ polynomials whereas with our method, we reduce the number of polynomials to 4. There is, however, a computational trade-off.
The difficulty in the former method is finding good candidate polynomials with low trace, whereas the largest difficulty in the authors' method is computing the support of the optimal distribution given by their linear programming problem.

%

\subsection{Method of Proof}
 It is well-known that $I(\mu)$ is concave which means that for all probability measures $\mu_1$ and $\mu_2$,
\begin{equation}\label{cvx}
 I(\mu_1)+I(\mu_2)-2I(\mu_1,\mu_2)\leq 0,   
\end{equation}
where $I(\mu_1,\mu_2):=\int \log|x_1-x_2|d\mu_1(x_1)d\mu_2(x_2)$ ~\cite[Lemma 2.1]{MR2730573}. It follows from Lemma~\ref{linlem} that 
\begin{equation}\label{energy}
I(\mu)\geq 0    
\end{equation}
is equivalent to 
the following family of linear inequalities 
\begin{equation}
2I(\mu,\mu_i)-I(\mu_i)= \int (2U_{\mu_i}(x)-I(\mu_i) )d\mu \geq 0    
\end{equation}
for every measure $\mu_i$ with non-negative energy $I(\mu_i)\geq 0$, where
\[
U_{\mu}(x)=\int \log|x-z|d\mu(z)
\]
is the logarithmic potential of $\mu$. This implies that the following optimization problem:
\[
\begin{array}{cc}
\min \int x d\mu(x)      \\
  \int \log|Q|d\mu(x)\geq 0  \text{ for every } Q\in A \\
  I(\mu)\geq 0,  
\end{array}
\]
is equivalent to the following linear programming problem:
\[
\begin{array}{cc}
\min \int x d\mu(x) \\
  \int \log|Q|d\mu(x)\geq 0  \text{ for every } Q\in A\\  
  \int (2U_{\mu_i}(x)-I(\mu_i) )d\mu \geq 0,
\end{array}
\]
for every measure $\mu_i$ with non-negative energy $I(\mu_i)\geq 0$. We denote the above optimization problem by the \textit{primal problem}.   Let 
$ \Lambda_{A}$ be the minimum of the primal problem. Suppose that $\mu_A$ is the probability measure that achieves this minimum, then
\[
\int x d\mu_A(x)=\Lambda_{A}.
\]
The associated dual problem is maximizing $\lambda$ such that
\[
x\geq \lambda +\sum_{Q_i\in A}\lambda_{Q} \log|Q(x)|+\sum_{\mu_i} \lambda_i (2U_{\mu_i}(x)-I(\mu_i))
\]
for $x\ge 0$,
where $\mu_i$ varies among probability measures. Let 
$ \lambda_{A}$ be the maximum of the dual problem. It follows that $\Lambda_A\geq \lambda_A$. As part of the proof of our main theorem, we show that strong duality holds which means $\Lambda_A= \lambda_A$.
By integrating the above inequality with respect to $d\mu_A$, it follows that the above inequality is equality on the support of $\mu_A$ and $\mu_i=\mu_A$. 
Then
\begin{equation}\label{main_eq}
x\geq \lambda_A +\sum_{Q\in A}\lambda_Q \log|Q(x)|+\lambda_0 (2U_{\mu_A}(x)-I(\mu_A))
\end{equation}
for $x\ge 0$,
\[
\int \log|Q(x)| d\mu_A=0
\]
for $Q\in A$ if $\lambda_Q\neq 0$, and 
\[
I(\mu_A)=0.
\]
In fact, the above conditions identify $\mu_A$ uniquely. It follows under some assumptions that $\mu_A$ is supported on the finite union of some intervals $\bigcup_{i=0}^l[a_{2i},a_{2i+1}] $ and has the following density function up to a scalar:
\[
\frac{|p(x)|\sqrt{|H(x)|}}{\prod_{\lambda_Q\neq 0}|Q(x)|},
\]
where $H(x)=\prod_{i=0}^l(x-a_{2i})(x-a_{2i+1})$ and $p(x)$ is a polynomial with degree $\sum_{\lambda_Q\neq 0} \deg(Q)-l-1$.

\subsection{Organization of the Paper}
In section~\ref{strong duality}, we formulate a convex optimization problem and prove that strong duality holds for this problem. In section~\ref{den_numer}, we give an expression for the optimal solution of the optimization problem in terms of polynomials. In section~\ref{zeroden}, we give a proof of Theorem~\ref{main1}.  In section~\ref{special_case}, we give proofs of Corollary~\ref{schur_bound} and Corollary~\ref{siegel_bound}. Finally in section~\ref{numeric}, we discuss our gradient descent algorithm and present our numerical results.

\section{Convex optimization and strong duality}\label{strong duality}
In this section, we give a proof of the following proposition.
\begin{proposition}\label{mainprop}
Let $A=\{Q_1,\dots,Q_n \}$ be any finite subset of polynomials. Define
\[
\Lambda_A=\inf \int xd\mu(x)
\]
subjected to $I(\mu)\geq 0$ and the following constraints:
\[
\int \log|Q_i(x)| d\mu(x)\geq 0
\]
for every $Q_i\in A$. Then there exists a unique probability measure $\mu_A$  and constants $\lambda_Q,\lambda_0\ge 0$ such that $\int xd\mu(x)=\Lambda_A$, $I(\mu)=0$,  
\[
x\geq \Lambda_A+\sum \lambda_Q\log|Q_i(x)|+\lambda_0U_{\mu_A}(x)
\]
for every non-negative $x$ with equality for every $x$ in the support of $\mu_A$, and 
\[
\int \log|Q_i(x)| d\mu_A(x)= 0
\]
if $\lambda_Q\neq 0$. 
\end{proposition}

The proof of the above proposition follows from several auxiliary lemmas. We begin by showing the support of $\mu_A$ is inside a fixed compact set independent of $A$. Our proof is based on the Principle of Domination for potentials that we cite below from~\cite[Theorem 3.2]{logpotentials}. 
\newline

\begin{theorem}~\cite[Theorem 3.2]{logpotentials}\label{dom}
    Let $\mu$ and $\nu$ be two positive finite Borel measures with compact support on $\mathbb{C}$ and suppose that the total mass of $\nu$ does not exceed that of $\mu$. Assume further that $\mu$ has finite logarithmic energy. If the inequality
    \[
    U_{\mu}(z)\geq U_{\nu}(z)+c 
    \]
 holds $\mu$-almost everywhere for some constant $c$, then it holds for all $z\in \mathbb{C}$.
\end{theorem}

First, we cite the following lemma from~\cite[Lemma 2.1]{MR2730573}. 
\begin{lemma}~\cite[Lemma 2.1]{MR2730573}\label{convexity}
    Suppose that $\mu_1,\dots, \mu_n$ are probability measures with finite logarithmic energy. Let $\sum_{i=1}^n c_i=1$ and $c_i> 0$. We have
    \[
    I\left(\sum_i c_i\mu_i\right)\geq \sum_{i} c_iI(\mu_i)
    \]
    with equality if and only if all $\mu_i$ are identical. 
\end{lemma}
Define
\[
I(\eta,\mu):=\int U_{\eta}(x)d\mu(x).
\]
\begin{lemma}\label{linlem}
    Suppose that $\mu$ is a probability measure with finite energy. If $I(\mu)\geq 0$, then 
    \[
    I(\eta,\mu)-\frac{I(\eta)}{2}\geq 0
    \]
    for every probability measure $\eta$ with finite energy. On the other hand, if for every $\eta$ with
    $I(\eta)= 0$,
    \[
    I(\eta,\mu)\geq 0,
    \]
   then $I(\mu)\geq 0.$
\end{lemma}
\begin{proof}
    Suppose that $I(\mu)\geq 0$. By  Lemma~\ref{convexity}, we have  
    \[
    I\left(\frac{\mu+\eta}{2}\right)\geq \frac{I(\mu)+I(\eta)}{2}.
    \]
The above implies
    
    \[
    I(\eta,\mu)-\frac{I(\eta)}{2}\geq \frac{I(\mu)}{2}\geq 0
    \]
    for every probability measure $\eta$ with finite energy.
    On the other hand, suppose that 
    \[
     I(\eta,\mu)-\frac{I(\eta)}{2}\geq 0
    \]
    for every $\eta$ with $I(\eta)= 0$.
    Suppose to the contrary that $I(\mu)<0$. Let $\sigma$ be the equilibrium measure of the interval $[0,5]$.
    We note that $I(\sigma)=\log(5/4)> 0.$ So there exists a parameter $0<\alpha<1$ such that 
    \[
    I(\alpha\mu+(1-\alpha)\sigma)=\alpha^2I(\mu)+2\alpha(1-\alpha)I(\mu,\sigma)+(1-\alpha)^2I(\sigma)=0.
    \]
    Take $\eta:=\alpha\mu+(1-\alpha)\sigma$. Note that $I(\eta)=0$ and 
    \begin{equation}\label{eqsubI}
          -\frac{\alpha}{2} I(\mu)-\frac{(1-\alpha)^2}{2\alpha}I(\sigma)=(1-\alpha)I(\mu,\sigma).
    \end{equation}

     By our assumption, we have
    \[
     I(\eta,\mu)=\alpha I(\mu)+(1-\alpha)I(\mu,\sigma)\geq 0.
    \]
    By substituting~\eqref{eqsubI}, we obtain  
    \[
    I(\eta,\mu)=\alpha I(\mu)+(1-\alpha)I(\mu,\sigma)=\frac{\alpha}{2} I(\mu)-\frac{(1-\alpha)^2}{2\alpha}I(\sigma)\geq 0
    \]
    This is a contradiction as $I(\mu)<0$ and $I(\sigma)>0$. So we must have $I(\mu)\geq 0$, concluding our lemma. 
\end{proof}

\begin{lemma}
    Suppose that $\mu$ and $\nu$ are probability measures with finite energy, and for some $\alpha>0$,
    \[
    U_{\mu}(z)\geq \alpha U_{\nu}(z) 
    \]
    for every $z\in \mathbb{C}$.
    Then
    \[
    I(\mu)\geq \alpha^2 I(\nu).
    \]
\end{lemma}
\begin{proof}
    We have
    \begin{align*}
     I(\mu)= \int U_{\mu}(z)d\mu(z)
        &\geq \alpha \int U_{\nu}(z)d\mu(z)
            \\
         &= \alpha \int\int \log|x-z| d\nu(x) d\mu(z)
            \\
         &= \alpha \int U_{\mu}(x)d\nu(x)
         \\
         & \geq \alpha^2 \int U_{\nu}(x)d\nu(x)= \alpha^2 I(\nu).
   \end{align*}
This completes the proof of our lemma.
\end{proof}

\begin{lemma}\label{compact}
Let $A=\{Q_1,\dots,Q_n \}$ be any finite subset of integer polynomials. Suppose that $\mu$ is a probability measure on $\mathbb{R}^+$ such that $I(\mu)\geq 0$ and 
\(
\int \log|Q_i(x)| d\mu(x)\geq 0
\)
for every $Q_i\in A$. There exists a probability measure $\mu'$ supported on $[0,18]$ such that  $I(\mu')\geq 0$,  
\(
\int \log|Q_i(x)| d\mu(x)\geq 0,
\)
and 
\[
\int x d\mu'(x)\leq \int xd\mu(x),
\]
where the above is equality only if $\mu$ is supported on  $[0,18]$.
\end{lemma}
\begin{proof}
 Our method is to truncate the support of $\mu$ with a parameter to reduced its expected value and at the same time to take a convex combination with a fixed equilibrium measure to satisfy the positivity conditions on the energy and log moments.
 We now explain this in detail. If $\int x d\mu(x)\geq 2$, let $\mu'=\mu_{[0,4]}$ be the equilibrium measure on the interval $[0,4]$.
 This is sufficient for the lemma since the potential is non-negative on the complex plane. Otherwise, $E:=\int x d\mu(x)< 2$. For $M>0$, let 
\[d\mu_M(x):=\begin{cases}
    \frac{d\mu(x)}{\mu([0,M])} &\text{ if $x<M$,}
    \\
    0 &\text{ otherwise,}
\end{cases}
\]
and 
\[
E_M:= \int x d\mu_M(x).
\]
Here $\mu_M$ is the conditional probability  of $\mu$ with the condition $x<M$. It follows that $E_M\leq E$ with equality if $\mu_M=\mu$.
In fact, by writing $E$ as a convex combination of conditional expectations on $x\leq M$ and $x>M$, we have
\[
E=\mu([0,M]) E_M+ \mu ([M, \infty])E_{M,\infty},
\]
where $E_{M,\infty}:=\frac{1}{\mu ([M, \infty])}\int_M^{\infty} x d\mu(x)$.
It is clear that $E_M\leq M$ and $E_{M,\infty} >M$ unless $\mu ([M, \infty])=0$.
For $N>4 $ and $M>N$, 
define
\[
\mu_{M,N,\alpha}:=\alpha\mu_{M}+(1-\alpha)\mu_{[0,N]},
\]
where $\mu_{[0,N]}$ is the equilibrium measure on the interval $[0,N]$. We have 
\begin{align*}
    \int x d\mu_{M,N,\alpha}(x)=\alpha E_M+ (1-\alpha)N/2. 
\end{align*}
Next, we estimate the potential of $\mu_{M,N,\alpha}$. 
Our goal is to find parameters $M,N, \alpha$ such that
\begin{equation}\label{1-ineq}
    \alpha E_M+ (1-\alpha)N/2 \leq E,
\end{equation}
and
    \begin{equation}\label{2-ineq}
      U_{\mu_{M,N,\alpha}}(x) \geq \delta U_{\mu}(x)     
    \end{equation}
    for every $x\in \mathbb{R}^+$ and some $1\geq \delta>0$. Inequality~\eqref{1-ineq} is equivalent to
    \[
    \alpha\geq \frac{N/2-E}{N/2-E_M}.
    \]
 We take $\alpha:=\frac{N/2-E}{N/2-E_M}$. So inequality~\eqref{1-ineq} is equality with this choice.  By~Theorem~\ref{dom}, 
    it is enough to check inequality~\eqref{2-ineq} for every $0\leq x<M$ as $0 \leq \delta\leq 1$. 
Suppose that $0\leq x<M$, we have 
\begin{align*}
    U_{\mu}(x)&= \int_0^M \log |x-z| d\mu(z)+\int_M^{\infty} \log |x-z| d\mu(z)
    \addtocounter{equation}{1}\tag{\theequation} \label{U-mu}\\
    &=\mu([0,M])U_{\mu_M}(x)+ \int_M^{\infty} \log |x-z| d\mu(z)
    \\
    & \leq \mu([0,M])U_{\mu_M}(x)+ \int_M^{\infty} \log |z| d\mu(z).
\end{align*}
Let $U_{\mu_{[0,N]}}(x)$ be the potential of $\mu_{[0,N]}$. We have 
\begin{equation}\label{U-mn}
       U_{\mu_{M,N,\alpha}}(x)= \alpha U_{\mu_M}(x)+(1-\alpha)U_{\mu_{[0,N]}}(x)\geq \alpha U_{\mu_M}(x)+ (1-\alpha)\log(N/4).
\end{equation}
By \eqref{U-mu} and \eqref{U-mn}, inequality~\eqref{2-ineq} follows if 
\begin{equation}\label{3-ineq}
\alpha U_{\mu_M}(x)+ (1-\alpha)\log(N/4)\geq \delta\left(\mu([0,M])U_{\mu_M}(x)+ \int_M^{\infty} \log |z| d\mu(z)\right).    
\end{equation}
We take $\delta:=\frac{\alpha}{\mu([0,M])}$. We show that $0 \leq \delta\leq 1$ which is equivalent to 
\[
\mu([0,M])\geq \alpha= \frac{N/2-E}{N/2-E_M}.
\]
This     is equivalent to 
\[
 N/2\leq \frac{1}{\mu([M,\infty])} \int_{M}^{\infty}x d\mu(x), 
\]
which follows from $N/2\leq M$. By this choice of $\delta$, inequality~\eqref{3-ineq}, is equivalent to 
\[
(1-\alpha)\log(N/4)\geq \delta \int_M^{\infty} \log |z| d\mu(z),
\]
Substituting the value of $\alpha$ in the above,  it remains to check 
\[
\frac{E-E_M}{N/2-E_M}\log(N/4)\geq \delta \int_M^{\infty} \log |z| d\mu(z).
\]
By substituting the values of $E$ and $E_M$, This is equivalent to 
\[
\int_{0}^{\infty} x d\mu(x)- \frac{1}{\mu([0,M])} \int_{0}^{M} x d\mu(x) \geq \int_M^{\infty}  \delta \frac{N/2-E_M}{\log(N/4)}\log |z| d\mu(z),
\]
or
\[
\frac{1}{\mu([M,\infty])}\int_{M}^{\infty} z-\delta \frac{N/2-E_M}{\log(N/4)}\log |z| d\mu(z) \geq \frac{1}{\mu([0,M])}\int_{0}^{M} x d\mu(x)=E_M.
\]
Note that $E_M\leq  E <2$. The above inequality follows if  
\[
 z-\delta \frac{N/2-E_M}{\log(N/4)}\log |z| \geq 2
\]
for every $z\geq M$. Let $N=4e$. Then, the above follows if
\[
z\geq 2e \log(z)+2.
\]
The above holds for every $z\geq 18$.
This concludes the proof of our lemma. 
\end{proof}

\begin{lemma}\label{lem2.5}
   Recall the definition of $\Lambda_A$ from Proposition~\ref{mainprop}. There exists a probability measure $\mu_A$ supported in $[0,18]$ satisfying all the above constraints and $ \int xd\mu_A(x) = \Lambda_A.$
\end{lemma}
\begin{proof}
    By definition of $\Lambda_A$ and Lemma~\ref{compact}, there exists a sequence of probability measures $\{\mu_n\}$ that are supported on the interval $[0,18]$ satisfying 
    the constraints
\[I(\mu_n)\ge 0,\]
\[
\int \log|Q(x)| d\mu_n(x)\geq 0
\]
for every $Q\in A$,
and
\[
\lim_{n\to \infty} \int xd\mu_n(x)=\Lambda_A.
\]
By compactness of the space of probability measures supported on $[0,18]$, there exists a sub-sequence $\mu_{a_n}$ that converges weakly to some probability distribution $\mu_0$.
Note that $\log |Q(z_1,\dots,z_m)|$ is an upper semi-continuous function. Hence, by monotone convergence theorem, we have 
\[
\int \log |Q(z_1,\dots,z_m)|d\mu_0(z_1)\dots d
\mu_0(z_m) \geq \limsup_{n\to \infty }\int \log |Q(z_1,\dots,z_m)|d\mu_{a_n}(z_1)\dots d
\mu_{a_n}(z_m)\geq 0.
\]
This completes the proof of our lemma. 
\end{proof}

\subsection{Strong duality}
In this section, we give a proof of Proposition~\ref{mainprop}. Recall the primal problem and its dual from section~\ref{prime} of the Introduction.
Let $\mathcal{P}$ be the space of Borel probability measures supported on $\mathbb{R}^+$ with finite energy. Let $\mathcal{P}^+\subset \mathcal{P}$ be the subspace of the probability measures with logarithmic energy greater or equal to zero on $\mathbb{R}^+$.
Let $A=\{Q_1,\dots,Q_n \}$ be a finite subset of integral polynomials.  By Lemma~\ref{lem2.5}, there exists a feasible solution $\mu_A \in \mathcal{P}^+$ with expectation $\Lambda_A$ to the primal problem. Next, we introduce the associated dual problem.  
\newline

Define  $\lambda_A$ to be the supremum over all $\lambda$ such that for every positive real $x$ and some fixed  elements $\mu_j\in \mathcal{P}^+$ and constants $c_i,d_j\geq 0$ then
\[
x\geq \lambda +\sum_{Q_i\in A}c_i \log|Q_i(x)|+\sum_{\mu_j} d_j (2U_{\mu_j}(x)-I(\mu_j)).
\]
By integrating the above inequality with respect to $\mu_A$ which satisfies the constraints of the primal problem, it follows that 
\[
\Lambda_A\geq \lambda_A.
\]
The above inequality is known as the weak duality, and if equality holds it is known as the strong duality. Proposition~\ref{mainprop} implies  that $\Lambda_A= \lambda_A$ and there exists a feasible optimal solution for the dual problem.
Finally, we give a proof of Proposition~\ref{mainprop} which is similar to the proof of strong duality from~\cite[section 5.3.2]{Boyd}.
\begin{proof}[Proof of Proposition~\ref{mainprop}] 
Let $\mathcal{P}_{1}^+\subset \mathcal{P}^+$ be the countable subset of smooth probability measures with density function $p(x)dx$ on finite union of intervals $\Sigma\subset \mathbb{R}^+$ and energy greater or equal to zero, where $p(x)\in \mathbb{Q}[x]$ up to a scalar and the end points of $\Sigma$ are rational numbers.  Note that $\mathcal{P}_{1}^+$  is dense inside $\mathcal{P}^+$  with respect to the weak topology.
We note that $\Lambda_A\leq 2$ since the equilibrium measure of the interval $[0,4]$ satisfies all the conditions of the primal problem and has expected value $2$.  We write
\[
\mathcal{P}_{1}^+:=\{\mu_i: i\in \mathbb{N}\}.
\]
For $N\geq 1$, we define $\psi_N:\mathcal{P}\to \mathbb{R}^{1+n+N}$ as follows:
\[
\psi_N(\mu):=\left(E_{\mu}(x),E_{\mu}(\log|Q_1|),\dots,E_{\mu}(\log|Q_n|),E_{\mu}(U_{\mu_1}(x)-I(\mu_1)/2),\dots, E_{\mu}(U_{\mu_N}(x)-I(\mu_N)/2)  \right) 
\]
where $E_\mu(f):=\int f(x)d\mu(x)$ is the expected value of  $f$ with respect to  $\mu$. Let 
\[\mathcal{A}_N:=\psi_N(\mathcal{P}) \subset  \mathbb{R}^{1+n+N}\]
be the image of $\psi_N$. For $\varepsilon>0$, let 
\[
\mathcal{B}_{\varepsilon,N}:=\{(x_0,x_1,\dots,x_{n+N}): x_0\leq \Lambda_A-\varepsilon,\text{ and }  x_i\geq 0 \text{ for }1\leq  i\leq n+N \}.
\]
Since $\mathcal{P}_{1}^+$ is dense inside $\mathcal{P}^+$, it follows from the definition of $\Lambda_A$ that for every $\varepsilon>0$, there exists a large enough $M$ such that 
\[
\mathcal{A}_N \cap \mathcal{B}_{\varepsilon,N} =\emptyset
\]
for every $N\geq M.$ Otherwise, there exists a sequence of probability measures $\mu_N\in \mathcal{P}$ for any $N\geq 0$ such that 
\[
E_{\mu_N}(x) \leq \Lambda_A-\varepsilon <2,
\]
\[
E_{\mu_N}(\log|Q_i|)\geq 0,
\]
for every $1\leq i\leq n$
and 
\[
E_{\mu_N}(U_{\mu_j}(x))\geq I(\mu_j)/2\geq 0
\]
for every $j\leq N$. Since $E_{\mu_N}(x) <2$,  $\{\mu_N\}$ is a tight family of probability measures on $\mathbb{R}^+$. By  Prokhorov's theorem, there exists a convergent subsequence $\{N_a\}$. Suppose that
\[
\mu_0=\lim_{a\to \infty} \mu_{N_a}.
\]
 It follows that 
\[
E_{\mu_0}(x) \leq \Lambda_A-\varepsilon,
\]
\[
E_{\mu_0}(\log|Q_i|)\geq 0,
\]
for every $1\leq i\leq n$,
and 
\[
E_{\mu_0}(U_{\mu_j}(x))\geq I(\mu_j)/2 \geq  0
\]
for every $\mu_j\in \mathcal{P}_{1}^+$. Next, we show that $I(\mu_0)\geq 0$, which gives a contradiction as $\mu_0$ is a feasible solution with expected value less than $\Lambda_A$. 
By Lemma~\ref{linlem}, it is enough to show that for every probability measure $\eta$ supported on $\mathbb{R}^+$ with $I(\eta)=0$, we have $I(\mu_0,\eta)\geq 0$. By density of $\mathcal{P}_{1}^+$ inside $\mathcal{P}^+$, there exists a sequence $\mu_{a_n}\in \mathcal{P}_{1}^+$ such that 
\[
\lim_{n\to \infty} \mu_{a_n}=\eta.
\]
This implies that 
\[
U_{\eta}(x)\geq \limsup_{n\to \infty} U_{\mu_{a_n}}(x)
\]
for every $x\in \mathbb{C}$. This implies that
\[
I(\eta,\mu_0)=E_{\mu_0}(U_{\eta}(x))\geq \limsup_{n\to \infty} E_{\mu_0}(U_{\mu_{a_n}}(x)) \geq 0.
\]
Hence, there exists a large enough $M$ such that 
\[
\mathcal{A}_N \cap \mathcal{B}_{\varepsilon,N} =\emptyset
\]
for every $N\geq M.$ 
Suppose that $N\geq M$. Note that $\mathcal{A}_N, \mathcal{B}_{\varepsilon,N}\subset \mathbb{R}^{1+n+N} $ are both convex regions with \[
\mathcal{A}_N \cap \mathcal{B}_{\varepsilon,N} =\emptyset.
\]
By the separating hyperplane theorem, it follows that there exists some coefficients $(c_0,\dots,c_n,d_1,\dots, d_N)$ such that
\begin{equation}\label{hyp}
  \sum_{i=0}^n c_i x_i +\sum_{j=1}^N d_j x_{j+n} \geq c_0E_{\mu}(x)+\sum_{i=1}^n c_iE_{\mu}(\log|Q_i|)+\sum_{j=1}^N d_jE_{\mu}(U_{\mu_j}(x)-I(\mu_j)/2)  
\end{equation}
for every $(x_0,x_1,\dots,x_{n+N})\in \mathcal{B}_{\varepsilon,N}$ and $\mu\in \mathcal{P}.$ Since $x_i\geq 0$ could be arbitrarily large for $i\geq 1$, so we must have $c_i,d_j\geq 0$ for every $i,j\geq 1$. Similarly, $x_0\leq\Lambda_A -\varepsilon $ can be arbitrarily small. So we must have $c_0\leq 0$. Next, we show that $c_0\neq 0.$ Let $\mu_{eq}$ be the
the equilibrium measure of $[0,18]$. We note that 
\[
E_{\mu_{eq}}(x)=9,
\]
\[
E_{\mu_{eq}}(\log|Q_i|)\geq \deg(Q_i)\log(9/2)
\]
and by inequality\eqref{cvx}, we have
\[
E_{\mu_{eq}}(U_{\mu_j}(x)-I(\mu_j)/2)\geq I(\mu_{eq})/2=\log(9/2)/2.
\]
By substituting $(x_0,x_1,\dots,x_{n+N})=(\Lambda_A-\varepsilon,0,\dots,0)$  and $\mu=\mu_{eq}$ in~\eqref{hyp} and the above identities, we have
\[
c_0(\Lambda_A-\varepsilon-9) \geq \sum_{i=1}^nc_i\deg(Q_i)\log(9/2)+\sum_{j=1}^N d_j\log(9/2)/2.
\]
Since $\Lambda_A\leq 2$, the above inequality implies that $c_0\neq 0$ and $c_0<0$. We divide both side of \eqref{hyp}, by $c_0$ and obtain
\[
\Lambda_A-\varepsilon \leq E_{\mu}(x)+\sum_{i=1}^n c_i'E_{\mu}(\log|Q_i|)+\sum_{j=1}^N d_j'E_{\mu}(U_{\mu_j}(x)-I(\mu_j)/2)
\]
for every $\mu\in \mathcal{P}$, where $c_i':=\frac{c_i}{c_0}\leq 0$ and $d_j':=\frac{d_j}{c_0}\leq 0$. Fix $x\in\mathbb{R}$ and let $\delta_x$ be the Dirac mass at $x$. Since there exists a sequence $\{\mu_k\}\subset \mathcal{P}$ such that $\lim_{k\to \infty} \mu_k=\delta_x$, it follows that 
\begin{equation}\label{mainineq}
    \Lambda_A-\varepsilon \leq  x+\sum_{i=1}^n c_i'\log|Q_i|(x)+\sum_{j=1}^N d_j'\left(U_{\mu_j}(x)-I(\mu_j)/2\right)
\end{equation}
for every $x\in \mathbb{R}^+$, and some coefficients $c_i',d_i'\leq 0$.
Let 
\[
\mu_0:=\frac{1}{d'}\sum_{j}d_j'\mu_j
\]
where $d'=\sum_{j}d_j'$. By linearity of the logarithmic potential and  Lemma~\ref{convexity}, we have 
\[
\sum_{j=1}^N d_j'\left(U_{\mu_j}(x)-I(\mu_j)/2\right) \leq d' \left(U_{\mu_0}(x)-I(\mu_0)/2  \right).
\]
By substituting the above in \eqref{mainineq}, we have 
\begin{equation}
    \Lambda_A-\varepsilon \leq x+\sum_{i=1}^n c_i'\log|Q_i|(x)+d' \left(U_{\mu_0}(x)-I(\mu_0)/2  \right)
\end{equation}
for every $x\in \mathbb{R}^+$ and some $c_i'\leq 0$ and $d' \leq 0$. By integrating the above inequality with respect to $\mu_{eq}$, we obtain 
\[
\Lambda_A-10\leq \Lambda_A-\varepsilon -9 \leq \sum_{i=1}^n c_i'\deg(Q_i)\log(9/2) +d'\log(9/2)/2.
\]
Define the compact subset $\Omega\subset \mathbb{R}^{n+1}$ as follows:
\[
\Omega:= \{(x_0,\dots,x_n): \Lambda_A-10\leq x_0 \log(9/2)/2+  \sum_{i=1}^{n} x_i \log(9/2)\deg(Q_i), x_i\leq 0 \text{ for every } i \}.
\]

Note that coefficients $c_1',\dots,c_n'$ and $d'$ depend on $\varepsilon > 0$ and $\Omega$ is independent of $\varepsilon.$ By definition of $\Omega$, we have $(d',c_1',\dots,c_n')\in \Omega$. Now we can take $\varepsilon_m=\frac{1}{m}$ and obtain a sequence of points $p_m:=(d_m',c_{1,m}',\dots,c_{n,m}')\in \Omega.$ Let $(D,C_1,\dots,C_n)$ be a limiting point of a sub-sequence. By~\eqref{mainineq}, we obtain
\begin{equation}\label{mainp}
    \Lambda_A-\frac{1}{m} \leq  x+\sum_{i=1}^n C_i\log|Q_i|(x)+D\left(U_{\mu_m}(x)-I(\mu_m)/2\right) 
\end{equation}
for some sequence of probability measures $\mu_m$. By Lemma~\ref{lem2.5}, there exists a feasible solution $\mu \in \mathcal{P}$ with optimal trace $\Lambda_A$ to the primal problem. We integrate the above inequality with respect to $\mu$ and obtain
\begin{equation}\label{eq}
    \Lambda_A-\frac{1}{m} \leq  \Lambda_A+\sum_{i=1}^n C_iE_{\mu}(\log|Q_i|(x))+DE_{\mu}\left(U_{\mu_m}(x)-I(\mu_m)/2\right)
\end{equation}
for every $m\geq 0$. Note that $E_{\mu}(\log|Q_i|(x))\geq 0$, $C_i\leq 0, D\leq 0$, and by ~\cite[Lemma 2.1]{MR2730573},
\[
E_{\mu}\left(U_{\mu_m}(x)-I(\mu_m)/2\right) = -I(\mu-\mu_m)/2+ I(\mu)/2,
\]
where 
\[
-I(\mu-\mu_m)=\lim_{R\to \infty}\int\int_{|z|<R}\left[ \int \frac{d\mu(x) -d\mu_m(x)}{|z-x|}\right]^2d\vol(z)\geq 0.
\]
By taking $m\to \infty$ in inequality~\eqref{eq}, we have
\begin{equation}\label{mainin}
0\leq \sum_{i=1}^n C_iE_{\mu}(\log|Q_i|(x))+D\left(I(\mu)/2+\limsup_{m\to \infty} -I(\mu-\mu_m)/2\right).    
\end{equation}
We note that the right-hand side is the sum of non-positive terms. It follows that \[
I(\mu)=0,
\]
\[
E_{\mu}(\log|Q_i|(x))=0
\]
if $C_i\neq 0$, and
\[
\limsup_{m\to \infty} -I(\mu-\mu_m)=0
\]
which implies $\lim_{m\to \infty}\mu_m=\mu$. This implies that 
\[
0=I(\mu)\geq \limsup_{m\to\infty}I(\mu_m).
\]
On the other hand, by our assumptions $I(\mu_m)\geq 0.$
So the following limit holds
\[
\lim_{m\to\infty}I(\mu_m)=0.
\]
This implies that 
\[
 \lim_{m\to \infty}U_{\mu_m}(x)=U_{\mu}(x)
\]
for every $x$ in the support of $\mu$.
From the above,~\eqref{mainp}, and \eqref{mainin}, it follows that 
\[
    \Lambda_A \leq  x+\sum_{i=1}^n C_i\log|Q_i|(x)+D\left(U_{\mu}(x)\right)
\]
with equality for every $x$ inside the support of $\mu$.
\end{proof}
Next, we show that the support of $\mu_A$ is a finite union of intervals. This proves part of  Theorem~\ref{main1}.
\begin{proposition}\label{fsupp}
    Let $\mu_A$ be the optimal solution to the primal problem in Proposition~\ref{mainprop} with support $\Sigma\subset [0,18]$ and $A$ being a finite subset of integer polynomials with only real roots. Then $\Sigma=\bigcup_{i=0}^l[a_{2i},a_{2i+1}] $ and every gap interval $(a_{2i+1},a_{2i+2})$ contains a root of some $Q\in A$.
\end{proposition}
\begin{proof}
    By Proposition~\ref{mainprop} and Lemma~\ref{lem2.5}, there exists a unique probability measure $\mu$ satisfying all the conditions of Theorem~\ref{main1}, but we need to check the support of $\mu$ is a finite union of intervals separated by the roots of $Q\in A$.
    Let $\Sigma$ be the support of $\mu$.
    Suppose that $a,b\in \Sigma$ and $Q(x)$ does not have any root inside the interval $[a,b]$ for every $Q\in A$. We claim that $[a,b]\subset \Sigma.$ Now suppose the contrary that there is an open interval gap inside $\Sigma\cap [a,b]$ say $(c,d)$ where $c$ and $d$ are maximal, i.e. $c,d\in \Sigma$ and $(c,d) \cap \Sigma=\emptyset$.
    Our argument is based on the following perturbation.
    Let $\delta>0$ be a small enough parameter that we choose later. Define 
    $\mu_1:=\mu|(c-\delta,c)$ and $\mu_2:=\mu|(d,d+\delta)$ to be the restriction of $\mu$ to two small intervals around $c$ and $d$. Let 
    \[
    0<m<\min(\mu_1((c-\delta,c)),\mu_2((d,d+\delta)))
    \]
    be smaller than the minimum of the two masses which by our assumption is strictly positive that we choose later.  Let $\mu_3$ be the equilibrium probability measure of the interval $[c-\delta,d]$.  Let $\mu'$ be the probability measure obtained by taking $m/2$ mass of $\mu_1$ and $\mu_2$ and replacing it with the $m\mu_3$, namely 
    \[
    \mu':=\mu-\lambda_1-\lambda_2+\lambda_3
    \]
    where $\lambda_1:=\frac{m}{2\mu_1((c-\delta,c))} \mu_1$, $\lambda_2:=\frac{m}{2\mu_2((d,d+\delta))} \mu_2$ and $\lambda_3:=m\mu_3$. 
    We have
    \[
    E_{\mu'}(x)-E_\mu(x)=-E_{\lambda_1}(x)-E_{\lambda_2}(x)+E_{\lambda_3}(x)< m \left(\frac{-c+\delta}{2}-d/2+\frac{d+c-\delta}{2} \right)=0
    \]

  Suppose that $y<c<d$. Since $\lambda_1$ is supported on $[c-\delta,c]$, we have 
\[
U_{\lambda_1}(y)\leq m \max(\log(c-y),\log(\delta))/2, 
\]
and similarly,
\[
U_{\lambda_2}(y)\leq m\log(d-y+\delta)/2=m(\log(d-y)/2+O(\delta)),
\]
where the implicit coefficient in $O(\delta)$ depends only on $d-c$. Hence, we have
    \begin{equation}\label{potinq}
    \begin{split}
        U_{\mu'}(y)- U_{\mu}(y)&=-U_{\lambda_1}(y)-U_{\lambda_2}(y)+U_{\lambda_3}(y)
        \\
        &\geq m\left(\log\left( \frac{c-y+d-y+2\sqrt{(c-y)(d-y)}}{4}\right)  -\log(c-y)/2-\log(d-y)/2 +O(\delta)\right) 
        \\
        &=m\left(\log \left(\frac{1}{2}+\frac{c-y+d-y}{4\sqrt{(c-y)(d-y)}}\right)+O(\delta)\right) > 0,
    \end{split}
    \end{equation}
 where we used the explicit value of the potential of the equilibrium measure $\mu_{3}$ and the fact that $\delta$ is small compared to $d-c$. 
Similarly, we have
\[
U_{\mu'}(y)-U_{\mu}(y)> 0
\]
 for $c<d<y.$ 
 Since every $Q\in A$ does not vanish inside $[a,b]$, this implies that 
    \[
    \sum_{Q(y)=0}U_{\mu'}(y)> \sum_{Q(y)=0}U_{\mu}(y)\geq 0
    \]
    for small enough $\delta$.
    Finally, we note that
    \begin{equation*}
    \begin{split}
        I(\mu')-I(\mu)= 2I(\mu,-\lambda_1-\lambda_2+\lambda_3)+I(-\lambda_1-\lambda_2+\lambda_3,-\lambda_1-\lambda_2+\lambda_3)
    \end{split}
    \end{equation*}
    By \eqref{potinq} and the fact that $\Sigma\cap (c,d)=\emptyset$,  we have
    \[
    I(\mu,-\lambda_1-\lambda_2+\lambda_3)=\int -U_{\lambda_1}(x)-U_{\lambda_2}(x)+U_{\lambda_3}(x) d\mu(x)=km > 0
    \]
    for some constant $k>0$.
    We also note that 
    \[
    I(-\lambda_1-\lambda_2+\lambda_3,-\lambda_1-\lambda_2+\lambda_3)=O(m^2).
    \]
    Therefore, by choosing $m$ small enough 
    \[
     I(\mu')>I(\mu).
    \]
We constructed $\mu'$ with expected value less than the optimal value $\Lambda_A$ and $\mu'$ satisfies all the constraints of the primal problem. This contradicts with the fact that the optimal expectation is $\Lambda_A$. This completes the proof of our main theorem.  
\end{proof}

\section{Explicit computation of logarithmic potential} \label{den_numer}

Fix a compact set $\Sigma:=\bigcup_{i=0}^l[a_{2i},a_{2i+1}]$ which is a finite union of closed intervals.  In this section, we construct measures with explicit potentials that we need to express the optimal solution to Theorem~\ref{main1}.
\subsection{Measures with locally constant potential}
 It is well-known that the equilibrium measure of $\Sigma$ has a constant logarithmic potential on $\Sigma$. Its density function is 
 \[
 f(x)dx=\frac{|p(x)|}{\pi\sqrt{|H(x)|}} dx
 \]
where $H(x)=\prod_{i=0}^{2l+1}(x-a_i)$ and $p(x)=\prod_{i=0}^{l-1}(x-\sigma_i)$ where $\sigma_i\in [a_{2i+1},a_{2i+2}]$ are the critical points of the Green's function of $\Sigma$  uniquely determined by the following equations: 
\begin{equation}\label{scmap}
\int_{a_{2i+1}}^{a_{2i+2}}\frac{p(x)}{\sqrt{|H(x)|}}dx=0
\end{equation}
for every  $0\leq i\leq l-1$~\cite{logpotentials}.
Generalizing this, we state the following lemma
which expresses the density function of measures with locally constant potentials. 
\begin{lemma}\label{localconstant}
    Suppose that $\mu$ is a measure with  density function 
 \[
 f(x)dx=\frac{|p(x)|}{\pi \sqrt{|H(x)|}} dx
 \]
where $H(x)=\prod_{i=0}^{2l+1}(x-a_i)$ and $p(x)=\prod_{i=0}^{l-1}(x-\sigma_i)$ where $\sigma_i\in [a_{2i+1},a_{2i+2}]$. Then the logarithmic potential of $\mu$ is locally constant and the difference of the value between consecutive intervals is
\[
U_{\mu}(a_{2i+2})-U_{\mu}(a_{2i+1})=\int_{a_{2i+1}}^{a_{2i+2}}\frac{p(x)}{\pi\sqrt{|H(x)|}}dx.
\]
\end{lemma}
\begin{proof}
It follows from the Schwarz-Christoffel map.
\end{proof}

\subsection{Measures with logarithmic potentials}
Suppose that $\alpha\notin \Sigma$ is given. The following proposition gives the measure with potential $\log(x-\alpha)$ on $\Sigma$. Let $\psi_\alpha(z):=\frac{1}{z-\alpha}$ which sends $\alpha$ to $\infty$ and let $\mu_\alpha$ be the equilibrium measure of $\Sigma_{\alpha}:=\psi_\alpha(\Sigma)$. Let $\nu_\alpha$ be the probability measure supported on $\Sigma$ which is the push-forward  of $\mu_\alpha$ with map $\psi_\alpha^{-1}(z):=\frac{1}{z}+\alpha$.   
\begin{lemma}\label{logpotlem}
Let $\nu_c$ be as above. We have
\[
U_{\nu_c}(x)=\log(x-\alpha)-U_{\mu_{\alpha}}(0)+C_{\Sigma_{\alpha}}.
\]
\end{lemma}
\begin{proof}
We have 
\begin{align*}
U_{\nu_\alpha}(x)=\int \log|x-y| d\nu_\alpha(y)&=\int \log\left|\alpha+\frac{1}{s}-\left(\alpha+\frac{1}{t}\right)\right| d\nu_\alpha\left(\alpha+\frac{1}{t}\right)
\\
&=\int \log\left|\frac{t-s}{st}\right| d\mu_{\alpha}(t) 
=\log\left(\frac{1}{s}\right)-U_{\mu_{\alpha}}(0)+C_{\Sigma_{\alpha}}.
\end{align*}
where we write $x=\alpha+\frac{1}{s}$ and $y=\alpha+\frac{1}{t}$ for some $t,s\in\psi(\Sigma)$. This completes the proof of our lemma. 
\end{proof}

\subsection{Measures with linear potentials}
\begin{lemma}\label{linear}
    Let $f(x)dx$ be a probability density on $\Sigma$ with a locally constant potential. Let $\eta$ be the signed measure with density $-xf(x)dx$ on $\Sigma$. Then the potential of $\eta$ is $x$ plus a locally constant function on $\Sigma$.
\end{lemma}
\begin{proof}
Suppose that $y\in \Sigma$. There exists some $i$ such that $y\in [a_i,b_i]$. Define 
\[
I_y:=[a_i,y]\subset \Sigma.
\]
By the assumption, we have 
\[
C_{i}=\int_{\Sigma} \log|x-t|f(x)dx
\]
for every $t\in I_y$ and some constant $C_i$.
We integrate the above identity over $I_y$ and obtain
\begin{equation}\label{int}
C_{i}(y-a_i)= \int_{\Sigma} \int_{a_i}^y \log|x-t|dt f(x)dx.    
\end{equation}
We integrate the inner integral and obtain 
\[
\int_{a_i}^y \log|x-t| dt= (y-x)\log|y-x|-(y-x) + (x-a_i)\log|x-a_i|-(x-a_i).
\]
We substitute the above in~\eqref{int} and obtain
\begin{align*}
C_{i}(y-a_i)&=\int_{\Sigma}\left[ (y-x)\log|y-x|-(y-x) + (x-a_i)\log|x-a_i|-(x-a_i) \right]f(x)dx.
\end{align*}
By simplifying the above identity, we deduce that
\[
y=\int_{\Sigma} -x\log|x-y|f(x)dx+c_i
\]
where 
\[
c_i:=a_i +\int_{\Sigma}x\log|x-a_i|f(x)dx
\]
is a constant that only depends on  $i$. 
\end{proof}
In order to make constant $c_i$ independent of $i$, we use Lemma~\ref{localconstant} to find a measure $\mu$ such that
\[
U_{\mu}(a_{2i+2})-U_{\mu}(a_{2i+1})=\beta(c_{i+1}-c_i)
\]
for some small enough $\beta>0$. By adding this measure to $\mu$, we construct a measure with linear potential.

\section{Proof of Theorem~\ref{main1}}\label{zeroden}
We give a proof of Theorem~\ref{main1} in this section. First, we show that the density function of the optimal probability measure to the primal problem  vanishes at the positive end points of the intervals in $\Sigma$.

\begin{proposition}\label{density_continuity}
Let $\mu$ be a solution to a primal problem. The density function of $\mu$ as a function on $(0,\infty)$ is continuous. In particular, it vanishes on the positive boundary points of $\Sigma$. 
\end{proposition}
\begin{proof}
By Proposition~\ref{mainprop}, there exists a unique measure $\mu_A$ which is supported on  $\Sigma\subset [0,18]$ such that
    \[
    \begin{array}{cc}
          \int x d\mu_A(x)=\lambda_{A},\\
         I(\mu_A)=0,          \\
                   x\geq \lambda_A +\sum_{Q\in A}\lambda_{Q} \log|Q(x)|+\lambda_0 U_{\mu_A}(x),\text{ and}\\
          \int \log|Q(x)| d \mu_A=0 \text{ if $\lambda_Q\neq 0$,}

    \end{array}
    \]
where in the third line the inequality holds for non-negative $x$ and equality holds for every $x\in \Sigma$ with some scalars $\lambda_Q\geq 0$ and $\lambda_0> 0$. By Proposition~\ref{fsupp},  $\Sigma=\bigcup_{i=0}^l[a_{2i},a_{2i+1}] $ is a finite union of intervals. Furthermore, every gap interval $(a_{2i+1},a_{2i+2})$ contains a root of some $Q\in A$. Since $\lambda_0\neq 0$, we have 
\begin{equation}\label{exppot}
U_{\mu_A}(x)= \frac{1}{\lambda_0}\left(x- \lambda_A -\sum_{Q\in A}\lambda_{Q} \log|Q(x)|\right)     
\end{equation}
for every  $x\in \Sigma$.
    Suppose to the contrary that the density function does not vanish at an end point $a_{i}$ for some $i$. We show this implies that the following inequality  
     \begin{equation}\label{maininq}
              x\geq \lambda_A +\sum_{Q\in A}\lambda_{Q} \log|Q(x)|+\lambda_0 U_{\mu_A}(x),
     \end{equation}
    fails to hold inside a neighborhood of $a_i$. This is a contradiction with Proposition~\ref{mainprop}.  Suppose that $a_i=a_{2k}$ for some $k$ where $a_i\neq 0$. Let $\delta>0$ be a small parameter that we choose later.
    By inequality~\eqref{maininq} and identity~\eqref{exppot}, we have 
    \[
    \frac{U(a_{2k}-\delta)-U(a_{2k})}{\delta}\leq \frac{1}{\lambda_0}\left(-1 +\sum_{Q\in A}\lambda_{Q} \frac{\log|Q(a_{2k})|-\log|Q(a_{2k}-\delta)|}{\delta} \right).
    \]
    This implies that 
    \begin{equation}\label{upperbd}
        \frac{U(a_{2k}-\delta)-U(a_{2k})}{\delta}\leq \frac{1}{\lambda_0}\left(-1 +\sum_{Q\in A}\lambda_{Q} \frac{Q'(a_{2k})}{Q(a_{2k})} \right)+O(\delta),
    \end{equation}
where the coefficients of $O$ depends only on $A$  and $\Sigma$. We note that if the density function of $\mu$ does not vanish at $a_{2k}$, then the density function goes to infinity near $a_{2k}$ with the following asymptotic behaviour 
\[
d\mu(x)\gg (x-a_{2k})^{-1/2}dx
\]
for $a_{2k}<x<a_{2k+1}$. This implies that
\[
\begin{split}
 \frac{U(a_{2k}-\delta)-U(a_{2k})}{\delta}&\gg \frac{1}{\delta}\int_{a_{2k}}^{a_{2k}+c}\left(\log(x-a_{2k}+\delta)-\log(x-a_{2k})\right)(x-a_{2k})^{-1/2}dx+O(1)
 \\
 &=\frac{1}{\delta}\int_{0}^{c}\log\left(1+\frac{\delta}{t}\right)t^{-1/2}dt+O(1)
 \\
  &=\delta^{-1/2}\int_{0}^{c/\delta}\log\left(1+\frac{1}{u}\right)u^{-1/2}du+O(1)
 \\
 &\gg \delta^{-1/2}+O(1). 
\end{split}
\]
where $t=\delta u$ and $a_{2k+1}-a_{2k}>c>0$.
The above is contradiction with \eqref{upperbd}. The case where $a_i = a_{2k+1}$ for some $k$ is analogous. This completes the proof of our theorem.
\end{proof}

Finally, we give a proof of our main theorem.
\begin{proof}[Proof of Theorem~\ref{main1}] Recall the notation in the proof of Proposition~\ref{density_continuity}.  Let \[
H(x)=\prod_{i=0}^l(x-a_{2i})(x-a_{2i+1}).
\]
By Lemma~\ref{linear} and Lemma~\ref{localconstant}, there exists a unique measure $\mu_{lin}$ with density function
\[
\frac{f_{lin}(x)}{\sqrt{|H(x)|}}
\]
for some polynomial $f_{lin}(x)\in \mathbb{R}[x]$ with $\deg(f)=l+1$
such that 
\(
U_{\mu_{lin}}(x)=x
\)
for every $x\in\Sigma$. Furthermore, by Lemma~\ref{logpotlem} and Lemma~\ref{localconstant}, for every $Q\in A$, there exists probability measure $\mu_{Q}$ with density function  
\[
\frac{f_Q(x)}{Q(x)\sqrt{|H(x)|}}
\] for some polynomial $f_{Q}(x)\in \mathbb{R}[x]$ with $\deg(f)=l+\deg(Q)$ such that
$U_{\mu_Q}(x)=\frac{\log|Q(x)|}{\deg(Q)}+C_Q$
for every $x\in \Sigma$ and some constant $C_Q$. By \eqref{exppot}, we have
\[
\mu_A=\frac{1}{\lambda_0}\mu_{lin}-\sum_{Q\in A}\frac{\lambda_Q\deg(Q)}{\lambda_0}\mu_Q+c\mu_{eq},
\]
where $c$ is a constant and $\mu_{eq}$ is the equilibrium measure of $\Sigma.$ This implies that the density function of $\mu_A$ is 
\[
\frac{|f(x)|}{\sqrt{|H(x)|}\prod_{\lambda_Q\neq 0}|Q(x)|},
\]
where $f(x)\in\mathbb{R}[x]$ with $\deg(f)=\sum_{\lambda_Q\neq 0} \deg(Q)+l+1$.  By Proposition~\ref{density_continuity}, $f(x)=H(x)p(x)$ for some polynomials $p(x)$ with $\deg(p)=\sum_{\lambda_Q\neq 0} \deg(Q)-l-1$
then $\mu_A$ has the following density function up to a scalar:
\[
\frac{|p(x)|\sqrt{|H(x)|}}{\prod_{\lambda_Q\neq 0}|Q(x)|}.
\]
This concludes the proof of Theorem~\ref{main1}.
\end{proof}

Given the explicit density function proven in Theorem ~\ref{main1}, we can actually prove differentiability of the potential of these pareto optimal measures.
\begin{lemma}
    \label{pot_diff}
    Let $\mu$ be a probability measure supported on $\Sigma:=\bigcup_{i=1}^n[a_{2i},a_{2i+1}]$ with density function $$f(x):=\frac{|(x-a_j)P(x)|}{|Q(x)\sqrt{|H(x)|}|},$$ where $a_j$ is any boundary point of $\Sigma,$  $H(x) = \prod_{i=1}^n (x-a_{2i})(x-a_{2i+1}),$ and $P(x),Q(x)$ are some real coefficients polynomials with roots outside $\Sigma$.  Then
    $\mu$ has a differentiable potential function at boundary point $a_j$.
\end{lemma}

\begin{proof}
 The argument is the same for each end point. For simplicity, we prove it at $a_0$. We know that $U_\mu$ is right-differentiable at $a_0$ by our explicit form for $U_\mu'(x)$ in $\Sigma$. So it suffices to show that
    \[
        \lim_{h\to 0^+}\int_\Sigma\frac{\log|a_0+h-x|+\log|a_0-h-x| - 2\log|a_0-x|}{h}f(x)dx=0.
    \]
    The left-hand side of this is equal to
    \[
        L:=\lim_{h\to 0^+}\int_\Sigma\frac{\log\left|1-\frac{h^2}{(a_0-x)^2}\right|}{h}f(x)dx.
    \]
    Doing the substitution $z = 1-\frac{h^2}{(a_0-x)^2}$ and noticing that $f(x) \le C\sqrt{|x-a_0|}$ for some positive constant $C>0$, we see that $|L|$ is at most a constant multiple of
    \[
        \lim_{h\to 0^+}\sqrt{h}\int_{-\infty}^1\frac{|\log|z||}{\sqrt{|1-z|^3}}dx=0.
    \]
    This completes the proof.
\end{proof}

\section{Proving Special Case Bounds}\label{special_case}
In this section, we use the integrals computed in the Appendix \ref{computation} to prove Corollaries \ref{schur_bound} and \ref{siegel_bound}.
To do this, we compute the measures $\mu_\emptyset$ and $\mu_{\{x\}}$ that satisfy the requirements of 
Theorem \ref{main1} explicitly. Let $A \in\{\emptyset, \{x\}\}$. We need the following:
\begin{align}
    \lambda_A &= \int xd\mu, \label{trace_eq}\\
    I({\mu_A}) &= 0 \label{energy_eq},\text{ and }\\
    \int \log xd\mu_A &= 0\label{log_eq}&\text{for }A=\{x\}.
\end{align}
To find this measure, we use \eqref{main_eq}, \eqref{trace}, \eqref{log_int}, and \eqref{energy_comp}. To use these formulas, we need to make a choice of $\alpha, \beta$ and $\gamma$. Notice that a measure of the form given in section \ref{computation} has density function 
\[\frac{\alpha + \frac{\beta\sqrt{ab}}{x}+\gamma(m-x)}{\pi\sqrt{(b-x)(x-a)}}.\]
Since this density must be positive, we have that $\alpha + \frac{\beta\sqrt{ab}}{x}+\gamma(m-x)\ge 0$ which rearranges to
$$\gamma x^2 - (\gamma m+\alpha) x - \beta\sqrt{ab}\ge 0.$$
This must hold for all $x\in [a,b]$. We make the unique choice of $\alpha, \beta,$ and $\gamma$ that makes the density 0 on the boundary, i.e. at $a$ and $b$, subject to the constraint $\beta = 1-\alpha$. The equality case of this quadratic inequality gives
$$x = \frac{\gamma m +\alpha\pm \sqrt{(\gamma m + \alpha)^2+4\gamma\beta\sqrt{ab}}}{2\gamma}.$$
Thus we must have that $\frac{\gamma m + \alpha}{2\gamma} = m$ and $\frac{\sqrt{(\gamma m +\alpha)^2+4\gamma \beta\sqrt{ab}}}{2\gamma} = \frac{b-a}{2}$.
\newline

The first equality gives us $\alpha = m\gamma$. The latter inequality gives us
$$(\gamma m + \alpha)^2 + 4\gamma\beta\sqrt{ab} = (b-a)^2\gamma^2.$$
Using $\alpha = \gamma m $ and $\beta = 1-\alpha$, we get
$$4\gamma(1-\gamma m)\sqrt{ab} = -4ab\gamma^2.$$
Dividing by $4\gamma\sqrt{ab}$ and rearranging, we finally have
$\gamma = \frac{1}{m-\sqrt{ab}}$.
Just as we use $m$ to represent the arithmetic mean of $a$ and $b$, we shall use $g$ to represent the geometric mean of $a$ and $b$ so that
\begin{equation}\label{coefficients}\alpha = \frac{m}{m-g}, \beta = \frac{-g}{m-g},\text{ and }\gamma = \frac{1}{m-g}.\end{equation}
This change of variables from $(a,b)$ to $(m, g)$ will prove useful in reducing the complexity of \eqref{trace}, \eqref{log_int}, and \eqref{energy_comp} to prove \eqref{trace_eq}, \eqref{energy_eq}, and \eqref{log_eq}.
As such, we will rewrite \eqref{trace}, \eqref{log_int}, and \eqref{energy_comp} in terms of only $m$ and $g$. Rewriting \eqref{trace}, we get
$$E:=\int_a^b xd\mu = \frac{m+g}{2}.$$
Rewriting \eqref{log_int}, we get
$$L:=\int_a^b\log xd\nu = \frac{m+g}{m-g}\log\frac{m+g}{2}-\frac{2g}{m-g}\log g-1.$$
Rewriting \eqref{energy_comp}, we get
$$I:=\frac{1}{2}\log\frac{m^2-g^2}{4}-\frac{2g^2}{(m-g)^2}\log\frac{m+g}{2g}+\frac{3g-m}{2(m-g)}.$$
The last equality we need to hold is \eqref{main_eq}.
Looking at the constant terms of this expression, equation \eqref{pot_comp} gives us that
$$\lambda_A=\frac{-(\alpha\log(c/2)-\beta\log\frac{a+b+2\sqrt{ab}}{b-a}-\gamma m)}{\gamma}$$
assuming $I=0$. Using the change of variables, we have
$$\lambda_A = m-\frac{m-g}{2}\log\frac{m-g}{2}-\frac{m+g}{2}\log\frac{m+g}{2}.$$
Lastly, we rewrite the equalities $L=0, I=0, \lambda_A=E$ in terms of $E$ and $g$
We see that $L=0$ is equivalent to
\begin{equation}\label{siegel_eq1}
E\log E - E =g\log g-g.
\end{equation}
Similarly, $\lambda_A=E$ is equivalent to
$$E\log E-E=-(E-g)\log(E-g) - g.$$
Consequently, this equation and \eqref{siegel_eq1} gives
\begin{equation}\label{siegel_eq2}
    g\log g = (g-E)\log(E-g).
\end{equation}
We now show that \eqref{siegel_eq1} and \eqref{siegel_eq2} imply $I=0$.
\begin{align*}
    I &=\frac{1}{2}\log\frac{m^2-g^2}{4}-\frac{2g^2}{(m-g)^2}\log\frac{m+g}{2g}+\frac{3g-m}{2(m-g)}\\
    &= \frac{1}{2}\log(E-g)+\left(\frac{1}{2}-\frac{2g^2}{(m-g)^2}\right)\log E+\frac{2g^2}{(m-g)^2}\log g+\frac{3g-m}{2(m-g)}\\
    &= \left(\frac{1}{2}-\frac{g}{m-g}\right)\log(E-g)+\left(\frac{1}{2}-\frac{2g^2}{(m-g)^2}\right)\log E+\frac{3g-m}{2(m-g)}
\end{align*}
by \eqref{siegel_eq2}.
We note that $\frac{1}{2}(1-\frac{2g}{m-g})$ is not 0 because $m>g$. Dividing the last line above by this quantity gives
$$\log(E-g) + \left(1+\frac{2g}{m-g}\right)\log E-1=0$$
which follows from \eqref{siegel_eq1} and \eqref{siegel_eq2}.
\newline

Therefore in the case $A=\{x\}$, if we have our requirements that $\lambda_A=E, I=0$, and $L=0$ then \eqref{siegel_eq1} and \eqref{siegel_eq2} hold simultaneously with no other required relations.
We claim that this system of equations
$\begin{cases}
    g\log g - g = E\log E - E\\
    g\log g = (g-E)\log(E-g)
\end{cases}$
has a unique solution.
\newline

\begin{proposition}\label{siegel_bound_prop}
    The system of equations
\[\begin{cases}
    x\log x - x = y\log y - y\\
    x\log x = (x-y)\log(y-x)
\end{cases}\]
has at most one real solution.
\end{proposition}
\begin{proof}
Note that the top equation is equivalent to $\left(\frac{x}{e}\right)^\frac{x}{e}=\left(\frac{y}{e}\right)^{\frac{y}{e}}$ and that from the second equation, we must have $y>x$.
Furthermore, $f(x)=x^x$ has a strictly increasing derivative. So $x^x=c$ has at most two solutions for a given $c$. Since $f'(x)=0$ has the unique solution $x=\frac{1}{e}$ and $x<y$, any solution to the system must satisfy $x<1<y$. This means that for a given $0< x < 1$, there is a unique $y>1$ such that the first equation holds. Call this number $y_1(x)$. Taking $\frac{d}{dx}$ of both sides of the top equation, we see that $\frac{dy_1}{dx} = \frac{\log x}{\log y_1}$. Since $0<x<1<y_1$, we thus have that $y_1(x)$ is decreasing.
\newline

Similarly, for a fixed $0<x<1$, the second equation has exactly one possible corresponding $y$. This is because we have $(y-x)\log(y-x)=-x\log x$ has exactly one solution due to $-x\log x$ being positive and the fact that $z\log z>0$ precisely when $z>1$. Furthermore, $z\log z$ has strictly positive derivative for $z>1$. So as in the first equation, for each $0<x<1$, there is a unique $y$ satisfying the second equation. Call this $y_2(x)$. Taking $\frac{d}{dx}$ of both sides of the bottom equation, we get $\frac{dy_2}{dx}=1-\frac{\log(x)+1}{\log(y_2-x)+1}$. We saw that $z=y-x>1$ and so $y > x+1$. Hence $\frac{dy_2}{dx}>-\log (x)>0$. Hence $y_2(x)$ is increasing.
Since $y_1(x)$ is decreasing and defined only on $x\in(0,1)$ and $y_2(x)$ is increasing on $(0,1)$, it follows that the curves intersect at most once.
\end{proof}

\begin{proof}[Corollary \ref{siegel_bound}]
From \eqref{siegel_bound_eq}, we have that $(1+\nu)\log(1+\nu^{-1})+\frac{\log\nu}{1+\nu}=1$ has a unique solution and Siegel derived the lowerbound $T:=e(1+\nu^{-1})^{-\nu}$. We wish to show that $T=\lambda_{\{x\}^+}$. Choose $T'=e(1+\nu^{-1})^{-(1+\nu)}$. Then $(T',T)$ is a solution to the system of equations in Proposition \ref{siegel_bound_prop} and hence the only solution. Indeed, the first equation in the system is a simple identity, and the second equation is equivalent to \eqref{siegel_bound_eq}. From our discussion above, it follows that $T=E=\lambda_{\{x\}^+}$.
\end{proof}
\begin{proof}[Corollary \ref{schur_bound}]
We can re-use the computation done for Corollary \ref{siegel_bound}, except now we require that $\lambda_\emptyset=E, I=0$, and $\beta=0$ instead of $L=0$. Since $\beta=0$, by \eqref{coefficients} we have that $g=0$. Since $g = \sqrt{ab}=0$ and $a<b$, it follows that $a=0$.
Note: In Section \ref{computation}, we assume $a>0$ to avoid $g(x)\equiv 0$; however all other computations are identical without this assumption. From this, we derive that $\alpha = 1$ and $\gamma = \frac{1}{m}=\frac{2}{b}.$ Now using \eqref{energy_comp}, we have $\log\left(\frac{b}{4}\right) = \frac{1}{2}$
since $I=0$. Hence $b = 4\sqrt e$. Thus $\lambda_\emptyset = E$ implies $\lambda_\emptyset = \frac{m}{2} = \sqrt e$ as desired.
\end{proof}

\section{Numerical results}\label{numeric}
In this section, we give numerical lower bounds on the Schur-Siegel-Smyth trace problem. First, we compute explicitly the optimal measure $\mu_A$ in terms of its support $\Sigma_A$. Next, we explain our gradient descent algorithm for approximating $\Sigma_A$.
\subsection{Expressing $\mu_A$}\label{Xdef}
Recall our notations from Theorem~\ref{main1}. Fix 
\(
A=\{Q_1,\dots,Q_n\}
\)
and let $\Sigma_A:=\bigcup_{i=0}^l[a_{2i},a_{2i+1}]$ be the support of the optimal measure $\mu_A$ in Theorem~\ref{main1}. Recall from the proof of Theorem~\ref{main1} in section~\ref{zeroden} that $\mu_{eq}$, $\mu_{lin}$, $\mu_Q$ are  the equilibrium measure, the measure with linear potential $U_{\mu_{lin}}(t)=t$, and the probability measure with potential $U_{\mu_{Q}}(t)=\frac{\log(|Q(t)|)}{\deg(Q)}+C_Q$ for some constant $C_Q$ supported on $\Sigma_A$ respectively. By Lemma~\ref{logpotlem}, 
the density function of  $\mu_Q$ is given by the following explicit expression:
\[
d\mu_Q(x)=\frac{1}{\deg(Q)}\sum_{Q(\alpha)=0}\frac{|P_{\alpha}(x-\alpha)|\sqrt{H(\alpha)}}{|x-\alpha|\pi\sqrt{H(x)}},
\]
where $P(x):=x^{\deg(P_{\alpha})}p_\alpha(1/x)$ is the polynomial that appears in the numerator of the equilibrium measure of the  image of $\Sigma$ under the map $z\to \frac{1}{z-\alpha}$. Suppose that we have exactly one root of $Q\in A$ inside every gap, then $\deg(p)=0$ in Theorem~\ref{main1} and up to a scalar $\mu_A$ has the following density function
\[
\frac{\sqrt{|H(x)|}}{\prod_{\lambda_Q\neq 0}|Q(x)|}.
\]
By the proof of Theorem~\ref{main1} in section~\ref{zeroden}, we have the following identity
\[
\begin{split}
X_{eq}\frac{|p_{eq}(x)|}{\pi\sqrt{H(x)}}+X_{lin}\frac{-x|p_{eq}(x)|+|P_{gap}(x)|}{\pi\sqrt{H(x)}}+\sum_{Q}\frac{X_Q}{\deg(Q)}\sum_{Q(\alpha)=0}\frac{|P_{\alpha}(x-\alpha)|\sqrt{H(\alpha)}}{|x-\alpha|\pi\sqrt{H(x)}}    \\=\frac{c|H(x)|}{\sqrt{H(x)}\prod_{Q\in \tilde{A}} |Q(x)|},
\end{split}
\]
where $X_{eq},X_{lin}>0, X_Q<0$ are some coefficients, and $c$ is normalized such that the right-hand side of the equation is a probability measure on $\Sigma$ and $\tilde{A}=\{Q\in A: X_{Q}\neq 0\}$. It follows that the  expressions inside the absolute values in the above have the same signs, so we may take the absolute values outside the sums. By comparing the top coefficient of both side, we have
\begin{equation}\label{lin}
    X_{lin}=c\pi
\end{equation}

Moreover, by computing residues at $x-\alpha$, we deduce that
\[
\frac{-X_Q}{\deg(Q)\pi}=\frac{c\sqrt{|H(\alpha)|}}{|Q'_\alpha(\alpha)|\prod_{\stackrel{Q\in \tilde{A}}{Q(\alpha)\neq 0}}|Q(\alpha)|}
\]
where $Q_\alpha\in\tilde{A}$ is such that $Q_\alpha(\alpha)=0$. Fix $Q\in \tilde{A}$ and let ${\alpha_1,\dots,\alpha_{\deg(Q)}}$ be roots of $Q$. The above identity implies that
\[
\frac{H(\alpha_1)}{\left(Q'_\alpha(\alpha_1)\prod_{\stackrel{Q\in \tilde{A}|Q(\alpha_1)|}{Q(\alpha_1)\neq 0}}|Q(\alpha_1)|\right)^2}=\dots=\frac{H(\alpha_{\deg(Q)})}{\left(Q'_{\alpha_{\deg(Q)}}(\alpha_{\deg(Q)})\prod_{\stackrel{Q\in \tilde{A}|Q(\alpha_{\deg(Q)})|}{Q(\alpha_{\deg(Q)})\neq 0}}|Q(\alpha_{\deg(Q)})|\right)^2}.
\]
We may write the above as
\begin{equation}\label{residue}
X_Q=-\frac{\deg(Q)c\pi\sqrt{|\text{Res}(H,Q)^{1/\deg(Q)}|}}{|\text{Disc}(Q)|^{1/\deg(Q)}|\prod_{Q_i\neq Q}|\text{Res}(Q,Q_i)|^{1/\deg(Q)}|}.
\end{equation}

\subsection{Gradient descent algorithm}
To compute the appropriate $\Sigma$ for a given set $A$, the authors implemented a gradient descent algorithm.
Given $\Sigma=\bigcup_i[a_{2i}, a_{2i+1}]$, we are able to compute $X_{lin}$, $X_{eq}$ and $X_Q$ for each $Q\in A$ as in section \ref{Xdef}.
We then compute $\mu_{eq, \Sigma}$, $\mu_{lin, \Sigma}$, and $\mu_{Q,\Sigma}$ for each $Q\in A$ as shown in section \ref{den_numer}. Define
$$\mu_\Sigma := X_{eq}\mu_{eq,\Sigma}+X_{lin,\Sigma}\mu_{lin,\Sigma}+\sum_{Q\in A}X_Q\mu_{Q,\Sigma}.$$
With this density function, we can compute all of the required quantities of a particular measure. So the high-level algorithm is to perform gradient descent on $\Sigma$ numerically (where $a_i$'s are the parameters) minimizing the sum of squares of the energy, the value of the density function at the $a_i$'s, and the values of $\int\log|Q(x)|d\mu_\Sigma$ for $Q\in A$, because these values are 0 from Theorem \ref{main1}.
\newline

In order to compute $d\mu_{eq,\Sigma}$, 
we use the equalities in \eqref{scmap}. To achieve these equalities, we can compute $\int_{a_{2i}}^{a_{2i+1}}\frac{x^j}{\sqrt{H(x)}}dx$ for each $i,j$ and solve a linear system to get the coefficients of the polynomial in the numerator. Similar systems can be solved to compute $d\mu_{lin, \Sigma}$ and $d\mu_{Q,\Sigma}$ since they can be derived from equilibrium measures. To save computation time, the authors combine the systems at each iteration so only one system needs to be solved per density function.
\newline

Lastly, to decrease the time taken by each iteration, the authors parallelized the computation of the gradient. They did this using multiprocessing. By using a different process for numerically computing the partial derivative of the objective function on each endpoint of $\Sigma$, we can compute the gradient much faster. For this, we use \texttt{Pool} from the \texttt{multiprocessing} module in Python.
\newline

One issue in the implementation of this algorithm is the numerical instability of integrating functions of the form $\frac{|p(x)|}{\sqrt{H(x)}}$ near the endpoints of the intervals. To fix this problem, when integrating on $[a_i,a_{i+1}]$, we perform the change of variables $u=\frac{a_{i+1}-a_i}{2}\cos(x)+\frac{a_{i}+a_{i+1}}{2}$.
This removes the factor of $\sqrt{(a_{i+1}-x)(x-a_i)}$ from the denominator of the integrand and hence removes the singularities. This allows for more stable numerical integration and consequently greatly increases the speed of the algorithm.
\newline

The above change of variables fixed the issue with the speed of each iteration in the gradient descent, but the descent itself was still quite slow. We found that instead of squaring the density of the measure at the boundary, we should square the density at the boundary divided by the density of the equilibrium measure at the boundary. This has given a great empirical improvement in the rate of the descent.

\subsection{Validity of Results}
Suppose there is a finite subset $A\subset\mathbb{Z}[x]$, a finite union of intervals $\Sigma=\bigcup_{i=0}^l[a_{2i},a_{2i+1}]$, a Borel probability measure $\mu$ with support $\Sigma$, and $\lambda>0,\lambda_0>0,$ and $ \lambda_Q>0$ for all $Q\in A$ such that $x = \lambda + \sum_{Q\in A} \lambda_Q \log|Q(x)| +\lambda_0\left(U_\mu(x)-\frac{I(\mu)}{2}\right)$ for all $x\in \Sigma$.
A priori, this alone is not sufficient in order to indicate whether $\lambda$ is a lower bound to the Schur-Siegel-Smyth trace problem. It is an additional requirement that
\begin{equation}\label{mu_ineq}
x \ge \lambda + \sum_{Q\in A} \lambda_Q \log|Q(x)| +\lambda_0\left(U_\mu(x)-\frac{I(\mu)}{2}\right)\text{ for all }x\in \mathbb{R}^+.
\end{equation}
We prove that equality on $\Sigma$ is indeed sufficient for special classes of $A$.
\begin{proposition}\label{eq_sufficient}
    Let $A, \Sigma,\lambda, \lambda_0,$ and $\lambda_Q$ for all $Q\in A$ be defined as above. Suppose $\Sigma\subset [0,18]$, and $\mu_{eq}$ is the equilibrium of $\Sigma$. Define $\mathcal R$ to be the set $\{\alpha\in\mathbb R: \exists Q\in A,Q(\alpha)=0\}$.
    Suppose that $\mathcal R\cap (a_{2i+1},a_{2(i+1)})=\{r_i\}$ for each $i=0,\dots, l-1$ and $\mathcal R\cap \big([0,a_0)\cup (a_{2l+1},\infty)\big)=\{0\}$. Moreover, suppose that for some  $\delta>0$ and any boundary point $a_i$ of $\Sigma$, we have
    \[
    \delta_{i}:=\lim_{x\to a_i,x\in \Sigma} \frac{d\mu(x)}{d\mu_{eq}(x)}\leq \delta, 
    \]
     for any $x\in (a_{2i+1},r_i)$,
    \[
    \sum_{Q\in A}\sum_{Q(\alpha)=0}\frac{\lambda_Q}{(x-\alpha)^2} + \lambda_0\int \frac{f_{2i+1}(y)}{(x-y)^2}dy>0
    \]
     where $f_{i}(y)$ is the density function of $\mu-\delta_{i}\mu_{eq}$, and similarly for any $x\in (r_i,a_{2(i+1)})$
    \[
    \sum_{Q\in A}\sum_{Q(\alpha)=0}\frac{\lambda_Q}{(x-\alpha)^2} + \lambda_0\int_{\Sigma} \frac{f_{2(i+1)}(y)}{(x-y)^2}dy>0.
    \]
     
    Then $\lambda-\delta \lambda_0 \log(18)$ is a lower bound on $\lambda_A$.
\end{proposition}
\begin{proof}
Suppose that $\mu_A$ is the optimal measure in Theorem~\ref{main1} with support  $\Sigma_A\subset[0,18]$, and 
\[
\int x d\mu_A(x)=\lambda_A. 
\]

    It suffices to prove 
    \[
    x \ge \lambda-\delta \lambda_0 \log(18) + \sum_{Q\in A} \lambda_Q \log|Q(x)| +\lambda_0\left(U_\mu(x)-\frac{I(\mu)}{2}\right)\text{ for all }x\in [0,18].
    \]
    Since by taking the average of both sides with respect to $d\mu_A(x)$, it follows that
    \[
    \lambda_A\geq \lambda-\delta \lambda_0 \log(18).
    \]
    Let 
    $$g(x) := x - \lambda+\delta \lambda_0 \log(18) - \sum_{Q\in A} \lambda_Q \log|Q(x)| -\lambda_0\left(U_{\mu}(x)-\frac{I(\mu)}{2}\right).$$ 
    It remains to show that $g(x)\ge 0$ for $0\le x\le 18$. By assumption, $g(x)= \delta \lambda_0 \log(18) > 0$ on $\Sigma$. Let $\{r_i\}=\mathcal R\cap (a_{2i+1},a_{2(i+1)})$. We show that $g(x)\geq 0$ for any $x\in (a_{2i+1},r_i)$. The argument for other intervals is similar. By our assumption, there exists $\delta_{2i+1}<\delta$ such that the density function of $\mu-\delta_{2i+1}\mu_{eq}$ vanishes at $a_{2i+1}$. Let 
     \[
    f(x) := x - \lambda- \sum_{Q\in A} \lambda_Q \log|Q(x)| -\lambda_0\left(U_{\mu-\delta_{2i+1}\mu_{eq}}(x)-\frac{I(\mu)}{2}\right).
    \]
Note that $f(x)=\delta_{2i+1} C_{\Sigma}\lambda_0 \geq 0$ on $\Sigma$, where $C_\Sigma>0$ is the logarithmic capacity of $\Sigma$  and
\[
g(x)-f(x)=\delta \lambda_0 \log(18)-\delta_{2i+1}\lambda_0U_{\mu_{eq}}(x)\geq 0
\]
for $0\le x\le 18$.
It is enough to show that $f(x)\geq 0$ for any $x\in (a_{2i+1},r_i).$  The key reason that $f(x)\ge 0$ for $x\in \mathbb R^+\setminus\Sigma$ is that
    $$f''(x) = \sum_{Q\in A}\sum_{Q(\alpha)=0}\frac{\lambda_Q}{(x-\alpha)^2} + \lambda_0\int_{\Sigma} \frac{f_{2i+1}(y)}{(x-y)^2}dy>0,$$
     where $f_{2i+1}(y)$ is the density function of $\mu-\delta_{2i+1}\mu_{eq}$, which is positive  by our assumption. We can differentiate under the integral by Liebnitz rule since the complement of the support of $f$ is open.
    Thus $f(x)$ has strictly increasing derivative outside $\Sigma$ on intervals of continuity.
    Furthermore, by Lemma \ref{pot_diff}, $f'(x)=0$ at $a_{2i+1}$. By our assumption on $\mathcal R$, this implies $g(x)\ge 0 $ for all $x\ge 0$.
\end{proof}

Thus to prove the corollaries in section \ref{grad_desc}, by Proposition \ref{eq_sufficient} and the explicit construction of measures with prescribed potential, it suffices to show that the resulting density functions we give are indeed positive.
The numerators of the density functions are polynomial combinations of absolute values of polynomials.
Thus we simply need to verify positivity of the numerators between the roots of these polynomials via polynomial inequalities. This can be done for each corollary given proving the bounds claimed in section \ref{grad_desc}.

\bibliographystyle{alpha}
\bibliography{main}

\appendix
\begin{section}{Ideal Measures on an Interval} \label{computation}
This section serve as a reference to potential and energy computations of the set of measures considered in this paper on $[a,b]$ where $0 < a<b$. Use $m:=\frac{a+b}{2}$ to denote the midpoint of the interval and $c := \frac{b-a}{2}$ to denote half the length of the interval.
Set $f(x) := \frac{1}{\pi\sqrt{(b-x)(x-a)}}$, $g(x) = \frac{\sqrt{ab}}{\pi x\sqrt{(b-x)(x-a)}}$ and $h(x) = \frac{m-x}{\pi\sqrt{(b-x)(x-a)}}$.

\subsection{Integrals}
This section simply verifies that $f$ and $g$ are density functions (corresponding to $\mu_{[a,b]}$ and $\nu_{[a,b]}$ respectively) and that $\int h d\eta=0$ where $\eta$ is the Lebesgue measure.
\begin{itemize}
    \item $\int fd\eta=1$.\newline
    This holds because $fd\eta$ is the density of the equilibrium measure on $[a,b]$.

    \item $\int gd\eta=1.$\newline
    This holds because $\nu_{[a,b]}$ is the pull-back measure of the equilibrium measure on $[b^{-1},a^{-1}]$ under the reciprocal map.
    \item $\int hd\eta=0.$\newline
    This holds because
    $\int_a^b h(x)dx = \int_{-c}^c\frac{-u}{\pi\sqrt{c^2-u^2}}du$ where $u=x-m$. The latter integrand is odd and therefore the integral is 0.
\end{itemize}

We have that if $\alpha,\beta,\gamma\in\mathbb R$, then $\tilde f:= \alpha f + \beta g + \gamma h$ has integral 1 so long as $\alpha+\beta = 1$. For $\tilde f$ to be a density function, it remains to consider only $\alpha, \beta,$ and $\gamma$ for which $\tilde f \ge 0$.

\subsection{Potentials}
Now we compute the potential of such $\tilde f$. To do so, we integrate $f,g$, and $h$ against logarithm. For this section, we assume $x,y\in[a,b]$.
\begin{itemize}
    \item $U_{\mu_{[a,b]}}(x) = \int_a^b \log|x-y|f(y)dy=\log(c/2)$ by the definition of the equilibrium measure.
    \item $U_{\nu_{[a,b]}}(x) = \int_a^b\log|x-y|g(y)dy = \log x - \log\frac{a+b+2\sqrt{ab}}{b-a}$ by Lemma \ref{logpotlem}.
    \item Lastly, $\int \log|x-y|h(y)dy = x - m$. The computation for $a=-c$ and $b=c$ is computed in Lemma \ref{linear}.
\end{itemize}
Therefore for a valid choice of $\alpha, \beta,$ and $\gamma$, we have that
\begin{equation}\label{pot_comp}U_{\tilde\mu}(x) = \int \log|x-y|\tilde f(y) dy = \gamma x + \beta\log x + \left(\alpha \log(c/2) - \beta\log\frac{a+b+2\sqrt{ab}}{b-a}-\gamma m\right)\end{equation}
where $\tilde \mu$ is the measure with density function $\tilde f$ for $a\le x \le b$.

\subsection{Expectations}\label{expectations}
We compute the expectation of $\tilde\mu$ by using the integrals of $xf(x),xg(x)$, and $xh(x)$. These integrals will be helpful in computing the energies.
\begin{itemize}
    \item $\int xf(x)dx=m$.\newline
    We can see this as follows:
    \begin{align*}
        \int_a^b xf(x)dx &= \int_a^b \frac{x}{\pi\sqrt{(b-x)(x-a)}}dx\\
        &= \int_{-c}^c\frac{u+m}{\pi\sqrt{c^2-u^2}}du
        &\text{where }u=x-m.\text{ So}\\
        &= \int_{-c}^c\frac{m}{\pi\sqrt{c^2-u^2}}du
        &\because \frac{u}{\sqrt{c^2-u^2}}\text{ is odd. Thus}\\
        &= m.
    \end{align*}
    The last equality holds, because we already computed that $\int_{-c}^c \frac{du}{\pi\sqrt{c^2-u^2}}=1.$

    \item $\int xg(x)dx = \sqrt{ab}.$\newline
    This one is simpler because $xg(x) = \sqrt{ab}f(x)$. Thus $\int_a^b xg(x) = \sqrt{ab}$ as desired.

    \item $\int xh(x)dx = \frac{-c^2}{2}.$\newline
    This computation is slightly more involved. First, we start as usual.
    $$\int_a^b xh(x)dx = \int_a^b \frac{mx-x^2}{\pi\sqrt{(b-x)(x-a)}}dx = -\int_{-c}^c \frac{u^2+um}{\pi\sqrt{c^2-u^2}}du$$ where $u=x-m$. Recall that $\frac{u}{\pi\sqrt{c^2-u^2}}$ is odd, so it suffices to show that $\int_{-c}^c\frac{u^2du}{\pi\sqrt{c^2-u^2}}=\frac{c^2}{2}$.
    Since $\frac{u^2}{\sqrt{c^2-u^2}}$ is even, we can write
    $$\int_{-c}^c\frac{u^2}{\pi\sqrt{c^2-u^2}}du = \frac{2}{\pi}\int_0^c\frac{u^2}{\sqrt{c^2-u^2}}du = \frac{2}{\pi}\left(\int_0^c\frac{c^2}{\sqrt{c^2-u^2}}du-\int_0^c\sqrt{c^2-u^2}du\right).$$
    It remains to notice that
    $$\int_0^c \frac{c^2}{\sqrt{c^2-u^2}}du = c^2\arcsin(u/c)|_{u=0}^c=\frac{1}{2}\pi c^2$$
    and that
    $$\int_0^c\sqrt{c^2-u^2}du = \frac{1}{4}\pi c^2$$ since the latter represents the area of a quarter circle of radius $c$.
    Putting everything together, we get
    $$\int_a^b xh(x)dx = -\int_{-c}^c\frac{u^2}{\pi\sqrt{c^2-u^2}}du = \frac{-2}{\pi}\left(\frac{1}{4}\pi c^2\right) = \frac{-c^2}{2}.$$
\end{itemize}

An interesting observation is that the expectation of $\mu_{[a,b]}$ is the arithmetic mean of $a$ and $b$ whereas that of $\nu_{[a,b]}$ is the geometric mean of $a$ and $b$. In particular, decreasing $\alpha$ gives a non-increasing expectation. Furthermore, since $\int xh(x)dx<0$, increasing $\gamma$ decreases expectation. Finally, the expectation of $\tilde\mu$ is
\begin{equation}\label{trace}\int_a^b x d\tilde\mu = \alpha m + \beta \sqrt{ab} -\frac{\gamma c^2}{2} = \sqrt{ab} + \alpha(m-\sqrt{ab}) - \frac{\gamma c^2}{2}.
\end{equation}

\subsection{Logarithmic Integrals}
We now compute when $\int \log |x| d\tilde\mu \ge 0$. These results are also useful in computing the energies.
\begin{itemize}
    \item $\int \log x d\mu_{[a,b]} = U_{\mu_{[a,b]}}(0) = \log\frac{a+b+2\sqrt{ab}}{4}.$\\
    This is the value using the notation from \cite{Smith} which is taken from Serre in \cite{MR2428512}.
    
    \item $\int \log x d\nu_{[a,b]}= U_{\nu_{[a,b]}}(0) = \log\frac{4ab}{a+b+2\sqrt{ab}}.$\\
    This can be derived from $U_{\mu_{[a,b]}}(0)$ using the $z\mapsto 1/z$ map.
    
    \item $\int_a^b \frac{(m-x)\log x}{\pi \sqrt{(b-x)(x-a)}}dx = \sqrt{ab}-m.$\\
    First, we see that
    \begin{equation}\label{h_zero_pot}\int_a^b\frac{\log|y-x|dx}{\pi\sqrt{(b-x)(x-a)}}=2\log\frac{\sqrt{a-y}+\sqrt{b-y}}{2}\end{equation}
    for $y\le a$.
    Integrating the left-hand side of \eqref{h_zero_pot} with respect to $y$ gives (up to a constant)
    $$\int_a^b\frac{(y-x)\log|y-x|-y}{\pi\sqrt{(b-x)(x-a)}}dx=y\left(2\log\left(\frac{\sqrt{a-y}+\sqrt{b-y}}{2}\right)-1\right)-\int_a^b\frac{x\log|y-x|dx}{\pi\sqrt{(b-x)(x-a)}}.$$
    Integrating the right-hand side of \eqref{h_zero_pot} with respect to $y$ gives (up to a constant)
    $$\sqrt{a-y}\sqrt{b-y}-(a+b-2y)\log\left(\frac{\sqrt{a-y}+\sqrt{b-y}}{2}\right).$$
    Therefore, we have that
    $$\int_a^b\frac{x\log |y-x|dx}{\pi\sqrt{(b-x)(x-a)}} = (a+b)\log\left(\frac{\sqrt{a-y}+\sqrt{b-y}}{2}\right)-\sqrt{a-y}\sqrt{b-y}-y+C.$$
    Finally,
    $$\int_a^b\frac{(m-x)\log|y-x|dx}{\pi\sqrt{(b-x)(x-a)}} =\sqrt{a-y}\sqrt{b-y}+y-m$$
    where the constant is $-m$ because we know that at $y=a$ the integral should be $-c$. In particular, when $y=0$, this is $\sqrt{ab}-m$.
\end{itemize}
Putting these together, we have that
\begin{equation}\label{log_int}
\int_a^b \log xd\tilde\mu = \alpha\log\frac{a+b+2\sqrt{ab}}{4} + \beta \log\frac{4ab}{a+b+2\sqrt{ab}} +\gamma \left(\sqrt{ab}-m\right).
\end{equation}

\subsection{Energies}
We wish to compute the energy of $\tilde \mu$. We do this by integrating the potentials of $\mu_{[a,b]}$ and $\nu_{[a,b]}$ and also by computing $\int\int \log|x-y|h(x)h(y)dxdy$.
\begin{itemize}
    \item $I({\mu_{[a,b]}}) = \log(c/2)$.\newline
    This holds because $\int_a^b U_{\mu_{[a,b]}}(x)d\mu_{[a,b]}(x) = \int_a^b \log(c/2)d\mu_{[a,b]}(x)=\log(c/2)$.
    \item $I({\nu_{[a,b]}}) = U_{\nu_{[a,b]}}(0)-\log\frac{a+2\sqrt{ab}+b}{b-a}.$\newline
    To compute this, we simply use the potentials we computed.
    \begin{align*}
        I({\nu_{[a,b]}}) &= \int_a^b U_{\nu_{[a,b]}}(x)d\nu_{[a,b]}\\
        &= \int_a^b \frac{(\log x - \log \frac{a+2\sqrt{ab} + b}{b-a})\sqrt{ab}}{\pi x \sqrt{(b-x)(x-a)}}dx\\
        &= \int_a^b \frac{\sqrt{ab}\log x}{\pi x\sqrt{(b-x)(x-a)}}dx - \log\frac{a+2\sqrt{ab}+b}{b-a}\int_a^b\frac{\sqrt{ab}}{\pi x\sqrt{(b-x)(x-a)}}dx\\
        &= U_{\nu_{[a,b]}}(0) - \log\frac{a+2\sqrt{ab}+b}{b-a}
    \end{align*}
    as desired.
    \item $\int\int \log|x-y|h(x)h(y)dxdy = \frac{-c^2}{2}.$\newline
    To do this, we recall that
    $$\int_a^b \log|x-y|h(y)dy = x-m.$$
    Since $\int_a^b h(x)dx = 0$, we just need to show that
    $$\int_a^b (x-m)h(x)dx = \int_a^b xh(x)dx = \frac{-c^2}{2}.$$
    This is shown in section \ref{expectations}.
\end{itemize}

Finally, we have that
\begin{align*}
I({\tilde\mu}) &= \int_a^b \gamma x + \beta\log x + \left(\alpha \log(c/2) - \beta\log\frac{a+2\sqrt{ab}+b}{b-a}-\gamma m\right)d\tilde\mu\\
&= \gamma \int_a^bxd\tilde\mu + \int_a^b \beta\log x(\alpha f(x)+\beta g(x)+\gamma h(x))dx + \left(\alpha \log(c/2) - \beta\log\frac{a+b
+2\sqrt{ab}}{b-a}-\gamma m\right)\\
&= \log(c/2) - 2\beta^2\log\frac{m+\sqrt{ab}}{2\sqrt{ab}}
-2\beta\gamma(m-\sqrt{ab}) - \frac{\gamma^2 c^2}{2}. \addtocounter{equation}{1}\tag{\theequation} \label{energy_comp}
\end{align*}
\end{section}

\begin{section}{Details on Computed Measures}\label{numerical_measures}
This appendix gives details of the measures suggested in Corollaries \ref{schur_bound} to \ref{our_bound}.
Firstly, the exact density functions of Corollaries \ref{schur_bound} and \ref{siegel_bound} are already given in the statement of these corollaries and these bounds are proven in Section \ref{special_case}.
\begin{figure}[ht]
     \centering
     \begin{subfigure}{1.\textwidth}
         \centering
        \includegraphics[width=\textwidth]{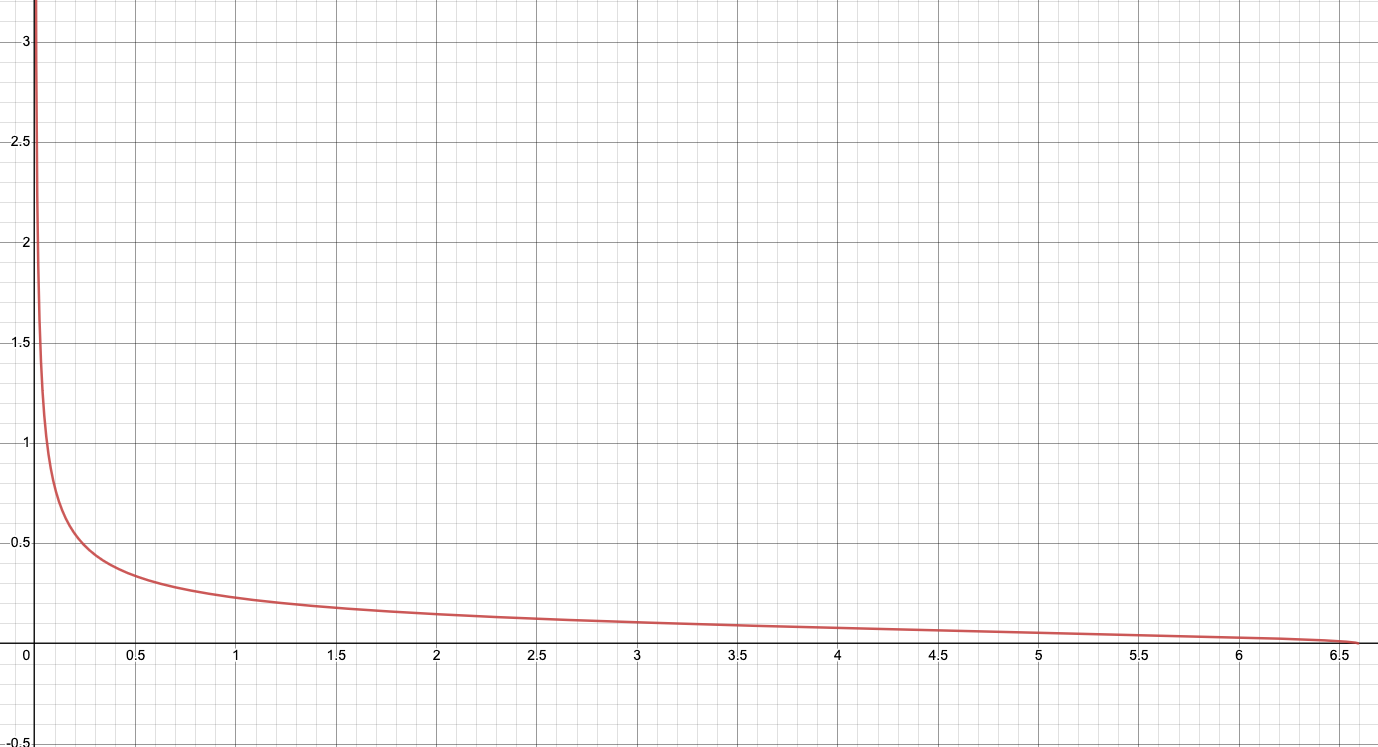}
         \caption{Density function of $\mu_{\emptyset}$}
        \label{schur_density}
     \end{subfigure}
     \vfill
     \begin{subfigure}{1.\textwidth}
         \centering         \includegraphics[width=\textwidth]{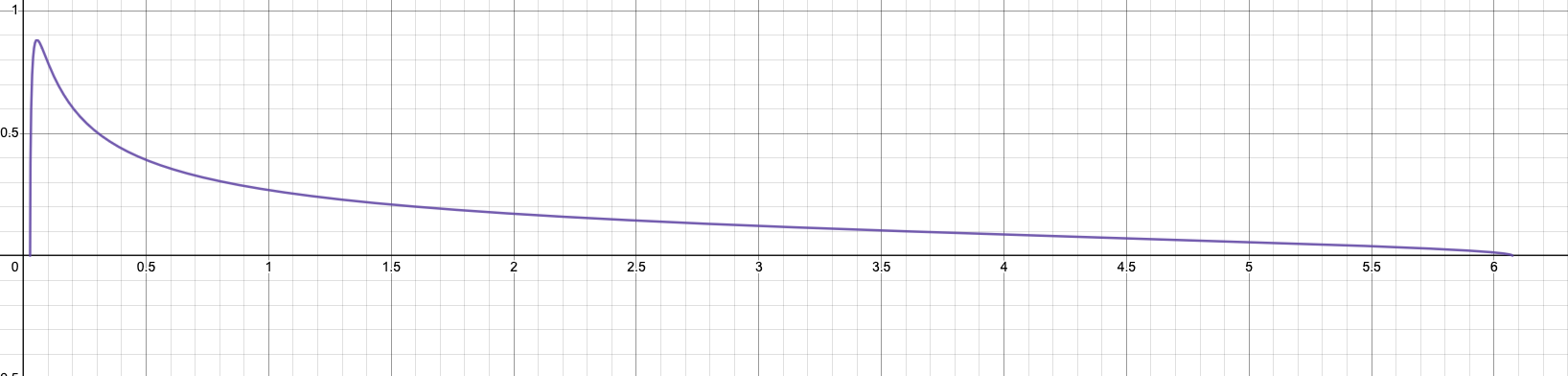}
                \caption{Density function of $\mu_{\{x\}}$}
        \label{siegel_density}
     \end{subfigure}
     \caption{Schur and Siegel Density Plots}
\end{figure}
The density functions for the cases of Schur's bound and Siegel's bound are given in Figures \ref{schur_density} and \ref{siegel_density}, respectively.
\newline

In this section, we give our numerical results for the choice of $A$ being $\{Q_0,Q_1,Q_2,Q_3,Q_4,Q_5,Q_6\}$ where we have $Q_0:=x,Q_1:=x-1, Q_2=x-2, Q_3:=x^2-3x+1,$ $Q_4:=x^{3}-5x^{2}+6x-1$, $Q_5:=x^4-7x^3+13x^2-7x+1$ and $Q_6:=x^4-7x^3+14x^2-8x+1$.
This gives a proof of Corollary~\ref{our_bound}. The proofs of Corollaries~\ref{2poly_bound}, \ref{3poly_bound}, \ref{4poly_bound}, and \ref{5poly_bound} are similar, and we just present the values of the supports of the measures which is enough to reproduce their proof. The supports for Corollaries \ref{2poly_bound} and \ref{3poly_bound} are in the statement of their corollaries.
The support for Corollary \ref{4poly_bound} is approximately
\begin{align*}
 &[0.0429949, 0.1859914]\cup [ 0.2104757, 0.3383178] \cup [0.4284029, 0.8045987]\cup[2208783, 1.5167041]\\
 &\qquad\cup[1.5933467, 2.4749011]\cup [2.7690575, 3.1791600]\cup [3.3160513, 5.4516406]
\end{align*}
and the support for Corollary \ref{5poly_bound} is approximately
\begin{align*}
    &[0.0471087, 0.1853294] \cup [0.2111668, 0.3366947]\cup [0.4302177, 0.8002212]\cup [1.2265017, 1.5147818]\\
    \cup &[1.5952680, 1.9682074] \cup[ 2.0322047, 2.4696444] \cup [2.7748974, 3.1757843] \cup [3.3195325, 5.3996085].
\end{align*}
\begin{proof}[Proof of Corollary~\ref{our_bound}]
Let  $\Sigma:=\bigcup_{k=0}^{15}[a_{2k},a_{2k+1}],$ where
$[a_0,\dots, a_{31}]:=$
\begin{align*}
 [&0.04865408503826852, 0.16838086459949675, 0.17793629526825372, 0.18381722534054817,\\
 &0.21253313678602553, 0.22282035304382763, 0.23253131212733313, 0.334984842461472,\\
 &0.43216747313335424, 0.5373228627969677, 0.550940510634415,0.6532454870830763,\\
 &0.670400493178262, 0.7958974217799853, 1.231938449021274, 1.512221832046952,\\
 &1.5978271338062664, 1.824427783903187, 1.8513776664372048, 1.9655369583386315,\\
&2.03493221184841, 2.192435667557317, 2.225841072618397, 2.463928802821324,\\
&2.7812624054781394, 3.1712064294805336, 3.324272345876191, 3.926025595570322,\\
&3.9869925948499505,4.354455187341252,4.426850255141888,5.35617801397137].
\end{align*}

By~\eqref{lin}, we let $X_{lin}$ to be
\[X_{lin}=0.89161282\approx c\pi. \] 
It is easy to check that the pairwise resultant of the above polynomials is $\pm 1$ (for distinct polynomials) except for $|\res(x, x-2)| = 2$. By \eqref{residue}, we choose to set the coefficients of the density functions to
\[
X_{Q_0}=-0.44801499\approx -\frac{c\pi|\res(Q_0,H)|^{1/2}}{2},
\]
\[
X_{Q_1}=-0.39003895\approx-c\pi|\res(Q_1,H)|^{1/2},
\]
\[
X_{Q_2}=-0.03655386\approx-\frac{c\pi|\res(Q_2,H)|^{1/2}}{2},
\]
\[
X_{Q_3}=-0.27820375\approx-\frac{2c\pi|\res(Q_3,H)|^{1/4}}{\sqrt{5}},
\]
\[
X_{Q_4}=-0.15505671\approx-\frac{3c\pi|\res(Q_4,H)|^{1/6}}{7^{2/3}},
\]
\[
X_{Q_5}=-0.05903767\approx-\frac{4c\pi|\res(Q_5,H)|^{1/8}}{725^{1/4}},
\]
\[
X_{Q_6}=-0.06244602\approx-\frac{4c\pi|\res(Q_6,H)|^{1/8}}{1125^{1/4}}.
\]
Finally, we let 
\[
X_{eq}=4.84960862
\]
and define
\[
\mu:=1.7(10)^{-5}\mu_{eq}+\frac{1.7(10)
^{-5}}{X_{eq}+X_{lin}\mu_{lin}(\Sigma)+\sum_{i=0}^3X_Q}\left(X_{eq}\mu_{eq}+\sum_{i=0}^3X_Q\mu_Q+X_{lin}\mu_{lin}\right).
\]

One can check that  for $x\in \Sigma$,
\[
x =  \lambda - \sum_{Q\in A} \lambda_Q \log|Q(x)| -\lambda_0\left(U_\mu(x)-\frac{I(\mu)}{2}\right)
\]
for every $x\in \Sigma$
where 
\[
\lambda \geq 1.802105.
\]
Next, we check the conditions of Proposition~\ref{eq_sufficient}. Let
 \[
    \delta_{i}:=\lim_{x\to a_i,x\in \Sigma} \frac{d\mu(x)}{d\mu_{eq}(x)}=\frac{P(a_i)}{P_{eq}(a_i)\prod_{Q\in A}Q(a_i)}.
    \]
One can check that $\delta_i\leq 3.7\times10^{-5}=:\delta$. Moreover, the inequality
    \[
    \sum_{Q\in A}\sum_{Q(\alpha)=0}\frac{\lambda_Q}{(x-\alpha)^2} + \lambda_0\int_{\Sigma} \frac{f_{2i+1}(y)}{(x-y)^2}dy\geq 0
    \]
     
    for $x\in (a_{2i+1},r_i)$  holds when
    \[
    \sum_{Q\in A}\sum_{Q(\alpha)=0}\frac{\lambda_Q}{(x-\alpha)^2} + \lambda_0\int_{\Sigma^+} \frac{f_{2i+1}(y)}{(x-y)^2}dy\geq -\lambda_0\int_{\Sigma^-} \frac{f_{2i+1}(y)}{(x-y)^2}dy,
    \]
   where  $\Sigma=\Sigma^+\cup\Sigma^{-}$ and
    \[
    \Sigma^+:=\{y\in \Sigma: f_{2i+1}(y)\geq 0 \}
    \]
    and 
    \[
    \Sigma^-:=\{y\in \Sigma:   f_{2i+1}(y)<0  \}.
    \]
    Note that $f_{2i+1}(y)dy=d(\mu-\delta_{2i+1}\mu_{eq})$, and we have
    \[
    \frac{\delta \lambda_0}{d(x,\Sigma^-)^2} > \delta_{2i+1}\lambda_0 \int_{\Sigma^-} \frac{d\mu_{eq}(y)}{(x-y)^2}dy>\left|\lambda_0\int_{\Sigma^-} \frac{f_{2i+1}(y)}{(x-y)^2}dy\right|,
    \]
    where $d(x,\Sigma^-):=\inf_{\sigma\in \Sigma^-}|x-\sigma|.$
    One can check numerically that $\Sigma^-$ lies inside a small neighbourhood of boundary points $a_{j}$ where $\delta_j< \delta_{2i+1}$. Indeed, $f_{2i+1}(a_{2i+1})=0$ and

    \[
    \frac{d}{dx}\frac{|P(x)|}{\prod_Q |Q(x)|}|_{x=a_{2i+1}}\leq \delta\cdot\frac{d}{dx} |P_{eq}(x)||_{x=a_{2i+1}},
    \]
    which implies $\Sigma^+$ contains a neighbourhood of $a_{2i+1}$. One can check numerically that 
    \[
\sum_{Q\in A}\sum_{Q(\alpha)=0}\frac{\lambda_Q}{(x-\alpha)^2} + \lambda_0\int_{\Sigma^+} \frac{f_{2i+1}(y)}{(x-y)^2}dy > \frac{\delta \lambda_0}{d(x,\Sigma^-)^2}
    \]
for every $x\in (a_{2i+1},r_i)$. By Proposition~\ref{eq_sufficient}, it follows that
\[
\lambda_A\geq 1.80203. 
\]
On the other hand, it is easy to check that $\mu$ satisfies all the constraints of Theorem~\ref{main1} and 
\[
\int xd\mu\leq1.80213.
\]
The above implies that
\[
\lambda_A\leq 1.80213.
\]
This completes the proof of  Corollary~\ref{our_bound}.
\end{proof}
\end{section}
\end{document}